\documentclass[11pt,a4paper,reqno]{amsart}

\usepackage[
    a4paper,
    left=2.8cm,
    right=2.8cm,
    top=2.5cm,
    bottom=2.5cm
]{geometry}

\usepackage{mathtools}
\usepackage{amssymb}
\usepackage{bm}
\usepackage{graphicx}
\usepackage{enumitem}

\usepackage[
    colorlinks=true,
    linkcolor=black,
    citecolor=black,
    urlcolor=blue
]{hyperref}

\usepackage{etoolbox}

\makeatletter
\patchcmd{\@setauthors}
  {\MakeUppercase{\authors}}
  {\authors\\[6pt]
   \normalfont\small
   \textsuperscript{a} School of Mathematics and Statistics and
   Chongqing Key Laboratory of Statistical Intelligent
   Computing and Monitoring, Chongqing Technology and
   Business University, Chongqing 400067, China.}
  {}{}
\makeatother

\theoremstyle{plain}
\newtheorem{theorem}{Theorem}[section]
\newtheorem{lemma}[theorem]{Lemma}

\theoremstyle{definition}

\theoremstyle{remark}

\numberwithin{equation}{section}

\newcommand{\R}{\mathbb{R}}

\newcommand{\Hs}{H^s(\R^3)}
\newcommand{\Hsr}{H_r^s(\R^3)}

\newcommand{\Dta}{\mathcal D_a^{t,2}(\R^3)}

\newcommand{\Dt}{\mathcal D^{t,2}(\R^3)}

\begin{document}

\title[Nontrivial solutions for a doubly fractional system]{Nontrivial Solutions for a Doubly Fractional\\Schr\"odinger--Bopp--Podolsky System in \texorpdfstring{$\R^3$}{R3}}

\author{Qi Liu\textsuperscript{a}}
\author{Lin Li\textsuperscript{a}}

\thanks{Corresponding author: Lin Li. Email: lilin420@gmail.com.}

\maketitle
\begin{abstract}
We study the following doubly fractional Schr\"odinger--Bopp--Podolsky system
\[
\begin{cases}
(-\Delta)^s u+Vu+q^2\varphi u=|u|^{p-2}u,\\
(-\Delta)^t\varphi+a^{2t}(-\Delta)^{2t}\varphi=u^2
\end{cases}
\quad\text{in }\mathbb R^3,
\]
where $a,V,q>0$, $0<s,t<1$, $4s+2t>3$ and $2<p<2_s^*$.
Here $2_s^*=6/(3-2s)$. If
$p>(4s+2t)/(s+t)$, we establish the existence of a
nontrivial solution for every $q>0$. For every
$2<p<2_s^*$, the same conclusion holds for
$0<q<q_0$, where $q_0$ is independent of $a$.
Moreover, in these parameter ranges, radial solutions can be
chosen so that, for every sequence $a_n\to0^+$,
a subsequence converges strongly in
$H_r^s(\mathbb R^3)\times\mathcal D_r^{t,2}(\mathbb R^3)$
to a nontrivial radial solution of the limiting
fractional Schr\"odinger--Poisson system.
\medskip


\medskip

\end{abstract}

\section{Introduction}\label{sec:intro}

The Schr\"odinger--Poisson system couples a nonlinear
Schr\"odinger equation with an electrostatic potential
generated by the density of the matter field. A typical
stationary model in $\mathbb R^3$ is given by
\begin{equation}\label{eq:classical-sp}
\begin{cases}
-\Delta u+Vu+q^2\varphi u=|u|^{p-2}u,
&x\in\mathbb R^3,\\
-\Delta\varphi=u^2,
&x\in\mathbb R^3,
\end{cases}
\end{equation}
where $u$ denotes the matter field, $\varphi$ is the
electrostatic potential, and $V,q>0$ are fixed
parameters. The positive sign of the electrostatic term corresponds
to a repulsive interaction. Schr\"odinger--Poisson equations with different
signs and physical normalizations also occur in optical and gravitational
models \cite{Paredes2020} and in semiconductor calculations
\cite{Hebal2021}. With the normalization in \eqref{eq:classical-sp}, the
Poisson equation determines the Coulomb potential
\[
\varphi_u(x)
=\frac{1}{4\pi}\int_{\mathbb R^3}
\frac{u^2(y)}{|x-y|}\,dy,
\]
whose kernel is singular at the origin.

To remove this singularity, Bopp and Podolsky
\cite{Bopp1940,Podolsky1942} introduced a higher-order
modification of classical electrodynamics. In the
electrostatic setting, the Poisson operator $-\Delta$
is replaced by
\[
L_a=-\Delta+a^2(-\Delta)^2,\qquad a>0.
\]
The corresponding equation is
\[
-\Delta\varphi+a^2(-\Delta)^2\varphi=u^2.
\]
Its fundamental solution is
\[
K_a(x)=\frac{1-e^{-|x|/a}}{4\pi|x|},
\qquad x\ne0,
\]
which extends continuously to the origin with
$K_a(0)=1/(4\pi a)$. Thus, the higher-order term
regularizes the Coulomb interaction at short distances.

Replacing the Poisson equation in
\eqref{eq:classical-sp} with the Bopp--Podolsky
equation leads to the classical
Schr\"odinger--Bopp--Podolsky system
\begin{equation}\label{eq:classical-sbp}
\begin{cases}
-\Delta u+Vu+q^2\varphi u=|u|^{p-2}u,
&x\in\mathbb R^3,\\
-\Delta\varphi+a^2(-\Delta)^2\varphi=u^2,
&x\in\mathbb R^3.
\end{cases}
\end{equation}
For the classical system, d'Avenia and Siciliano
\cite[Theorems~1.1--1.3]{DAveniaSiciliano} obtained a nontrivial solution
for $2<p<6$ and small charge, and for $3<p<6$ and arbitrary nonzero charge.
They also studied the radial limit as $a\to0^+$. Their normalization puts
$4\pi u^2$ on the right of the field equation, a constant that can be
absorbed into the coupling. Subsequent work includes critical nonlinearities
\cite{ChenTang,YangChenLiu,LiuChen}, multiplicity and nodal solutions
\cite{FigueiredoSiciliano,BahrouniMissaoui}, and nonconstant potentials
\cite{ChenLiRadulescuTang,JiangLiChen}. Other boundary and geometric
settings include the Neumann problem \cite{AfonsoSiciliano} and the Proca
system on compact manifolds \cite{Hebey2019,DAveniaGhimenti,Talluri}.

A different direction of development concerns the
introduction of fractional operators into the
Schr\"odinger--Poisson system. A typical fractional
Schr\"odinger--Poisson system has the form
\begin{equation}\label{eq:frac-sp-lit}
\begin{cases}
(-\Delta)^s u+\lambda\phi u=g(u),
&x\in\mathbb R^3,\\
(-\Delta)^t\phi=\lambda u^2,
&x\in\mathbb R^3,
\end{cases}
\end{equation}
where $0<s,t<1$. For general subcritical nonlinearities satisfying
Berestycki--Lions assumptions, Zhang, do \'O and Squassina
\cite[Theorem~1.1]{ZdOS} obtained positive radial solutions for small
$\lambda>0$ when $4s+2t\ge3$. As $\lambda\to0$, a subsequence tends
to $(U,0)$, where $U$ solves the uncoupled scalar equation. For equal
orders $s=t>3/4$, under suitable assumptions on the potential
and for a pure-power nonlinearity, Teng
\cite{Teng2019} obtained a ground state for $3<p<2_s^*$ in our
exponent convention. Critical growth, symmetry, decaying potentials and
prescribed mass have been studied in
\cite{Teng2016,YunAnYe,ZhangHe2023,HeMengSquassina}. In a critical
semiclassical Poisson problem with $2s+2t>3$, Yang, Yu and Zhao
\cite{YangYuZhao} used the lower exponent $(4s+2t)/(s+t)$.

The classical field operator has the factorization
$(-\Delta)(I+a^2(-\Delta))$. Replacing $-\Delta$ in its two factors
by $(-\Delta)^t$ gives $(-\Delta)^t(I+a^{2t}(-\Delta)^t)$.
Its inverse has an additional factor $(1+(a|\xi|)^{2t})^{-1}$ relative
to the Riesz potential. It suppresses high frequencies while retaining the
Poisson interaction at low frequencies. This leads to two questions.
For which powers does the resulting system have a nontrivial solution at
every fixed coupling? When the length parameter tends to zero with the
charge fixed, can one recover a nontrivial solution with a nonzero
electrostatic field? We consider the system
\begin{equation}\label{eq:sys}
\begin{cases}
(-\Delta)^s u+Vu+q^2\varphi u
=|u|^{p-2}u, &x\in\mathbb R^3,\\
(-\Delta)^t\varphi
+a^{2t}(-\Delta)^{2t}\varphi
=u^2, &x\in\mathbb R^3.
\end{cases}
\end{equation}
Here $0<s,t<1$, $a,V,q>0$, $4s+2t>3$, and
$2<p<2_s^*=6/(3-2s)$. These are the standing assumptions below.

For fixed $a>0$, we consider the completion of
$C_c^\infty(\mathbb R^3)$ with respect to the norm
\[
\|\varphi\|_a^2
=
\int_{\mathbb R^3}
\left(
|(-\Delta)^{t/2}\varphi|^2
+a^{2t}|(-\Delta)^t\varphi|^2
\right)\,dx.
\]
We denote this space by $\mathcal D_a^{t,2}(\mathbb R^3)$.
We refer to Section 2 for further properties of
this space.
For fixed $a,V,q>0$, we say that a pair
$(u,\varphi)\in H^s(\mathbb R^3)\times
\mathcal D_a^{t,2}(\mathbb R^3)$ is a weak
solution of system \eqref{eq:sys} if
\[
\begin{aligned}
&\int_{\mathbb R^3}
\left[
(-\Delta)^{s/2}u(-\Delta)^{s/2}v
+Vuv+q^2\varphi uv
\right]\,dx\\
&\qquad =
\int_{\mathbb R^3}|u|^{p-2}uv\,dx,
\quad \forall v\in H^s(\mathbb R^3),
\end{aligned}
\]
and
\[
\begin{aligned}
&\int_{\mathbb R^3}
\left[
(-\Delta)^{t/2}\varphi(-\Delta)^{t/2}\xi
+a^{2t}(-\Delta)^t\varphi(-\Delta)^t\xi
\right]\,dx\\
&\qquad =
\int_{\mathbb R^3}u^2\xi\,dx,
\quad
\forall \xi\in\mathcal D_a^{t,2}(\mathbb R^3).
\end{aligned}
\]
A solution $(u,\varphi)$ is called nontrivial if
$u\not\equiv0$.

Set
\[
p_*(s,t)=\frac{4s+2t}{s+t}.
\]
The first two results give existence at arbitrary coupling throughout
$p_*(s,t)<p<2_s^*$. We state the two ranges separately because the
estimate for the Palais--Smale sequence changes at $p=4$.

\begin{theorem}\label{thm:1.1}
If $s>\frac34$ and $p\in[4,2_s^*)$, then, for every
$q>0$, system \eqref{eq:sys} admits a nontrivial solution.
\end{theorem}

\begin{theorem}\label{thm:1.2}
If $p_*(s,t)<p<\min\{4,2_s^*\}$, then for every $q>0$
system \eqref{eq:sys} admits a nontrivial solution.
\end{theorem}

\begin{theorem}\label{thm:1.3}
For every $2<p<2_s^*$, there exists $q_0=q_0(s,t,p,V)>0$ such that,
for all $a>0$ and $q\in(0,q_0)$, system \eqref{eq:sys} admits a
nontrivial solution.
\end{theorem}

The three constructions can also be performed in the radial
subspace. The principle of symmetric criticality then yields
radial weak solutions of the full system. In the equal-order
case, the threshold satisfies $p_*(s,s)=3$, coinciding with
the lower exponent in the Schrödinger--Poisson problem
discussed above. The borderline case $p=p_*(s,t)$ is covered
by Theorem~\ref{thm:1.3} for sufficiently small $q>0$.

We next investigate the asymptotic behavior of
these radial solutions as $a\to0^+$. Formally
setting $a=0$ in system \eqref{eq:sys}, we obtain
the fractional Schr\"odinger--Poisson system
\begin{equation}\label{eq:sp}
\begin{cases}
(-\Delta)^s u+Vu+q^2\varphi u=|u|^{p-2}u,\\
(-\Delta)^t\varphi=u^2,
\end{cases}
\qquad x\in\mathbb R^3.
\end{equation}

\begin{theorem}\label{thm:1.4}
Fix $s,t,p,V,q$ as in one of Theorems~\ref{thm:1.1}--\ref{thm:1.3}.
For every $a>0$ one can choose a nontrivial radial solution
$(u_a,\varphi_a)$ of \eqref{eq:sys} with
$\sup_{a>0}\|u_a\|_{H^s}<\infty$. For each sequence $a_n\to0^+$
there are a subsequence, still denoted by $a_n$, and a nontrivial radial
solution $(u_0,\varphi_0)$ of \eqref{eq:sp} such that
\[
u_{a_n}\to u_0\text{ in }H^s_r,\qquad
\varphi_{a_n}\to\varphi_0\text{ in }\mathcal D^{t,2}_r,\qquad
a_n^t\|(-\Delta)^t\varphi_{a_n}\|_2\to0 .
\]
The same subsequential conclusion holds for any sequence of nontrivial
radial solutions at $a_n$ with $\sup_n\|u_{a_n}\|_{H^s}<\infty$.
\end{theorem}

Our approach is mainly based on variational methods.
Nevertheless, the doubly fractional structure of the
system gives rise to several analytical difficulties.
The main points are highlighted below.

(i) The first difficulty concerns the fundamental
solution of the electrostatic operator
$L_{a,t}=(-\Delta)^t+a^{2t}(-\Delta)^{2t}$.
Unlike the classical Bopp--Podolsky operator, whose
fundamental solution admits an explicit formula, the
fractional operator does not generally yield a
comparably explicit expression. We therefore establish
the positivity of its fundamental solution and the
uniform pointwise estimate
$0<K_{t,a}(x)\leq C_t|x|^{-(3-2t)}$
by means of complete monotonicity and Bernstein's
theorem. This estimate yields the Hardy--Littlewood--
Sobolev bounds needed to control the electrostatic
potential and the associated energy functional.

(ii) The different fractional orders $s$ and $t$
lead to distinct Sobolev structures for the matter
field and the electrostatic potential. In particular,
the condition $4s+2t>3$ ensures that the charge
density $u^2$ belongs to the natural source space
of the electrostatic equation. Another difficulty
arises in establishing the boundedness of
Palais--Smale sequences. Since the electrostatic
energy $T_a(u)=\int_{\mathbb R^3}\varphi_u u^2\,dx$
is homogeneous of degree four, the standard
energy-derivative estimate provides the required
bound when $p\geq4$, but fails when $2<p<4$,
because the electrostatic term appears with a
negative coefficient. To overcome this difficulty
for arbitrary charge, we introduce the stretching
family
$u_\tau(x)=\tau^{s+t}u(\tau x)$ and set
$d=4s+2t-3$. The Fourier multiplier gives, for
$\tau\geq1$,
\[
T_a(u_\tau)=\tau^d T_{a\tau}(u)
\leq \tau^d T_a(u),
\qquad
\|u_\tau\|_p^p
=\tau^{p(s+t)-3}\|u\|_p^p.
\]
The fractional kinetic term also scales as $\tau^d$.
Thus, the negative nonlinear term grows faster
precisely when
$p>p_*(s,t)=(4s+2t)/(s+t)$.
This comparison yields a common path of
bounded energy in the arbitrary-charge range.
Applying the abstract bounded Palais--Smale theorem
of Guida and Rolando
\cite[Theorem~1.1]{GuidaRolando}, together with their
quantitative estimate
\cite[Remark~2.2]{GuidaRolando}, we obtain bounds
for the selected subquartic Palais--Smale sequences
independently of $a$. For the small-charge result,
we instead truncate the quartic interaction and use
the estimate $T_a(u)\leq T_0(u)$, where $T_0$
denotes the Riesz interaction, to choose the
truncation and the charge independently of $a$.

(iii) In $H^s(\mathbb R^3)$, the lack of compactness
caused by translations makes the existence problem
more delicate. In the small-charge regime, the
norm-dependent truncation of the electrostatic
energy yields bounded Palais--Smale sequences.
To overcome the remaining lack of compactness,
we analyze these sequences by a splitting argument
adapted to the fractional Sobolev space and the
nonlocal electrostatic energy. A finite-profile
decomposition allows us to recover a nontrivial
critical point of the original functional.

(iv) Finally, the limit $a\to0^+$ presents a
further difficulty, since the higher-order term
in the electrostatic operator vanishes and the
associated spaces $\mathcal D_a^{t,2}$ vary with
$a$. Establishing uniform estimates for the
electrostatic potentials and the radial solution
families is therefore essential. The uniform control
of the sources $u_a^2$, together with the compact
embeddings of the radial fractional Sobolev space
$H^s_r$ and the corresponding dual embeddings,
allows us to pass to the limit in both equations.
Consequently, along a subsequence, the radial
solutions converge strongly to a nontrivial solution
of the limiting fractional Schrödinger--Poisson
system. The same compactness conclusion holds
for arbitrary uniformly bounded sequences of
nontrivial radial solutions.

The paper is organized as follows.
Section~\ref{sec:prelim} constructs the electrostatic
kernel and the reduced functional.
Section~\ref{sec:existence} proves the existence
theorems in the full space.
Section~\ref{sec:radial} treats the radial
construction and its limit as $a\to0^+$.
Appendix~\ref{app:splitting} proves the splitting
lemma used for the translation-invariant problem.

\subsection{Notations}

\begin{itemize}
\item For $r>1$, we denote by $\|\cdot\|_r$ the
norm in $L^r(\mathbb R^3)$, and by
$r'=r/(r-1)$ its conjugate exponent. We also
write $B_R$ for the ball of radius $R$ centred
at the origin.

\item Throughout the paper, $C,C_1,C_2,\ldots$
denote positive constants whose values may
change from line to line. Unless otherwise
specified, constants may depend on $s,t,a,V,p,q$ and the cutoff
functions. Uniformity in $a$ or $q$ is indicated where needed.

\item We denote the radial subspaces of
$H^s(\mathbb R^3)$,
$\mathcal D_a^{t,2}(\mathbb R^3)$ and
$\mathcal D^{t,2}(\mathbb R^3)$ by
$H^s_r(\mathbb R^3)$,
$\mathcal D_{a,r}^{t,2}(\mathbb R^3)$
and
$\mathcal D_r^{t,2}(\mathbb R^3)$,
respectively.
\item Unless otherwise indicated, all integrals
are taken over $\mathbb R^3$.
\end{itemize}

\section{Preliminaries}\label{sec:prelim}

\subsection{Fractional Sobolev spaces}

All function spaces for the weak system are real. We use complex Fourier
transforms to write the corresponding real Hilbert norms and pairings.
We fix the $L^2$-unitary Fourier transform
\[
\hat f(\xi)=(2\pi)^{-3/2}\int_{\R^3}e^{-ix\cdot\xi}f(x)\,dx,
\qquad
f(x)=(2\pi)^{-3/2}\int_{\R^3}e^{ix\cdot\xi}\hat f(\xi)\,d\xi,
\]
initially for Schwartz functions and then by the usual extension
to $L^2(\R^3)$. For $s\in(0,1)$, the fractional Laplacian
$(-\Delta)^s$ is defined as the Fourier multiplier with
symbol $|\xi|^{2s}$, namely,
\[
\widehat{(-\Delta)^s u}(\xi)
=|\xi|^{2s}\hat u(\xi).
\]
The fractional Sobolev space is defined as
\[
H^s(\mathbb R^3)
=\left\{
u\in L^2(\mathbb R^3):
\int_{\mathbb R^3}|\xi|^{2s}|\widehat u(\xi)|^2\,d\xi<\infty
\right\},
\]
endowed with the inner product
\[
(u,v)_{H^s}
=\int_{\mathbb R^3}
\left(
(-\Delta)^{s/2}u\,(-\Delta)^{s/2}v+Vuv
\right)\,dx,
\qquad V>0,
\]
and the corresponding norm
\[
\|u\|_{H^s}^2
=\int_{\mathbb R^3}
\left(
|(-\Delta)^{s/2}u|^2+V|u|^2
\right)\,dx.
\]
Set
\[
2_s^*=\frac{6}{3-2s},
\qquad
2_t^*=\frac{6}{3-2t}.
\]
We shall also use the homogeneous fractional Sobolev space
\[
\Dt = \left\{ u \in L^{2_t^*}(\R^3) : \int_{\R^3} |\xi|^{2t}|\hat u(\xi)|^2\,d\xi < \infty \right\},
\]
with inner product
\[
(u,v)_{\Dt} = \int_{\R^3} (-\Delta)^{t/2}u\,(-\Delta)^{t/2}v\,dx
= \int_{\R^3} |\xi|^{2t}\hat u(\xi)\overline{\hat v(\xi)}\,d\xi,
\]
and with norm
\[
\|u\|^2_{\Dt} = \int_{\R^3} |\xi|^{2t}|\hat u(\xi)|^2\,d\xi .
\]

The following embedding properties will be used
throughout the paper, particularly in establishing
the integrability of the nonlinear terms and the
compactness of radial sequences.

\begin{lemma}[Lions \cite{Lions82}]\label{lem:2.1}
Let $s\in(0,1)$. Then $H^s(\mathbb R^3)$ is continuously
embedded into $L^r(\mathbb R^3)$ for $2\leq r\leq 2_s^*$
and compactly embedded into $L^r_{\mathrm{loc}}(\mathbb R^3)$
for $1\leq r<2_s^*$. Moreover,
$H^s_r(\mathbb R^3)$ is compactly embedded into
$L^r(\mathbb R^3)$ for every $r\in(2,2_s^*)$.
\end{lemma}

\subsection{The electrostatic operator and its fundamental solution}\label{subsec:kernel}

For $a>0$ and $t\in(0,1)$, define
\[
L_{a,t}=(-\Delta)^t+a^{2t}(-\Delta)^{2t}.
\]
The associated energy space is $\Dta$, the completion
of $C_c^\infty(\R^3)$ with respect to the inner product
\[
(\varphi,\psi)_{\Dta}
=\int_{\R^3}
\left(
(-\Delta)^{t/2}\varphi\,(-\Delta)^{t/2}\psi
+a^{2t}(-\Delta)^t\varphi\,(-\Delta)^t\psi
\right)\,dx,
\]
and the corresponding norm
\[
\|\varphi\|_{\Dta}^2
=\int_{\R^3}
\left(
|(-\Delta)^{t/2}\varphi|^2
+a^{2t}|(-\Delta)^t\varphi|^2
\right)\,dx.
\]
For $\varphi\in C_c^\infty(\R^3)$, we have
\[
\|\varphi\|_{\Dt}\leq\|\varphi\|_{\Dta}.
\]
By density, this yields the continuous embedding
\[
\Dta\hookrightarrow\Dt.
\]
Furthermore, the fractional Sobolev inequality
\cite{CotsiolisTavoularis,DPV} gives
\[
\|u\|_{2_t^*}\leq S_t\|u\|_{\Dt},
\qquad u\in\Dt,
\]
for some constant $S_t>0$. Consequently,
\[
\Dta\hookrightarrow\Dt\hookrightarrow
L^{2_t^*}(\R^3).
\]
 We take the translation-invariant Green kernel with the inverse
Fourier multiplier displayed below. In particular,
\[
L_{a,t}K_{t,a}=\delta_0
\quad\text{in }\mathcal S'(\R^3).
\]
Under the Fourier convention fixed above, we have
\[
\widehat{K_{t,a}}(\xi)
=\frac{(2\pi)^{-3/2}}
{|\xi|^{2t}+a^{2t}|\xi|^{4t}}
=\frac{(2\pi)^{-3/2}}
{|\xi|^{2t}(1+a^{2t}|\xi|^{2t})}.
\]
Using the identity
\[
\frac{1}{|\xi|^{2t}(1+a^{2t}|\xi|^{2t})}
=
\frac{1}{|\xi|^{2t}}
-\frac{a^{2t}}{1+a^{2t}|\xi|^{2t}},
\]
we obtain, in $\mathcal S'(\R^3)$,
\begin{equation}\label{eq:Kdecomp}
K_{t,a}=R_t-G_{t,a}.
\end{equation}
Here
\[
\widehat{R_t}(\xi)
=(2\pi)^{-3/2}|\xi|^{-2t},
\qquad
\widehat{G_{t,a}}(\xi)
=(2\pi)^{-3/2}
\frac{a^{2t}}{1+a^{2t}|\xi|^{2t}}.
\]
Thus $R_t$ is the Riesz kernel of $(-\Delta)^t$.
The decomposition \eqref{eq:Kdecomp} will allow us
to establish the positivity and the Riesz-type
bound of $K_{t,a}$, which are essential for the
estimates below.

\begin{lemma}\label{lem:2.2}
The fundamental solution $K_{t,a}$ admits a
representative satisfying
\begin{equation}\label{eq:Rieszbound}
0<K_{t,a}(x)\leq c_t|x|^{-(3-2t)}
\quad\text{for a.e. }x\in\R^3\setminus\{0\},
\end{equation}
where
\[
c_t=
\frac{\Gamma\bigl((3-2t)/2\bigr)}
{2^{2t}\pi^{3/2}\Gamma(t)}.
\]
Consequently, by the Hardy--Littlewood--Sobolev
inequality \cite{Lieb},
\begin{equation}\label{eq:HLS}
\left|
\iint_{\R^3\times\R^3}
K_{t,a}(x-y)f(x)g(y)\,dx\,dy
\right|
\leq
C_t\|f\|_{(2_t^*)'}\|g\|_{(2_t^*)'}
\end{equation}
for all $f,g\in L^{(2_t^*)'}(\R^3)$,
where $C_t$ is independent of $a>0$.
\end{lemma}

\begin{proof}
Setting $\lambda=|\xi|^2$, let
\[
h_{t,a}(\lambda)=\frac{1}{\lambda^t(1+a^{2t}\lambda^t)},
\qquad
g_{t,a}(\lambda)=\frac{a^{2t}}{1+a^{2t}\lambda^t}.
\]
 The map $\lambda\mapsto\lambda^t$ is a Bernstein function.
Composing it with $z\mapsto(1+a^{2t}z)^{-1}$ gives a completely
monotone function. Multiplication by $\lambda^{-t}$ preserves complete
monotonicity. Bernstein's representation theorem
\cite[Chapters~1 and~3]{SchillingSongVondracek} therefore supplies
positive measures $\mu_{t,a}$ and $\nu_{t,a}$ such that
\[
h_{t,a}(\lambda)=\int_0^\infty e^{-\lambda\tau}\,d\mu_{t,a}(\tau),
\qquad
g_{t,a}(\lambda)=\int_0^\infty e^{-\lambda\tau}\,d\nu_{t,a}(\tau).
\]
Both transforms tend to zero as $\lambda\to\infty$, so neither
measure has an atom at $\tau=0$.
Writing $H_\tau(x)=(4\pi\tau)^{-3/2}e^{-|x|^2/(4\tau)}$, whose Fourier transform is $(2\pi)^{-3/2}e^{-\tau|\xi|^2}$, we obtain
\[
K_{t,a}(x)=\int_0^\infty H_\tau(x)\,d\mu_{t,a}(\tau),
\qquad
G_{t,a}(x)=\int_0^\infty H_\tau(x)\,d\nu_{t,a}(\tau)
\]
in the sense of distributions. Hence $K_{t,a}\geq0$ and $G_{t,a}\geq0$.

Moreover, by the monotone convergence theorem and Tonelli's theorem,
\[
\|G_{t,a}\|_{L^1}=\nu_{t,a}((0,\infty))
=\lim_{\lambda\to0^+}g_{t,a}(\lambda)=a^{2t}.
\]
Thus $G_{t,a}\in L^1(\R^3)$. Since $R_t\in L^1_{\mathrm{loc}}(\R^3)$, \eqref{eq:Kdecomp} implies $K_{t,a}\in L^1_{\mathrm{loc}}(\R^3)$. Therefore,
\[
K_{t,a}(x)\geq0,\qquad G_{t,a}(x)\geq0
\quad\text{a.e. in }\R^3.
\]
The heat-kernel representation and $\mu_{t,a}\neq0$ further give $K_{t,a}(x)>0$ almost everywhere. Consequently, by \eqref{eq:Kdecomp},
\[
0<K_{t,a}(x)\leq R_t(x)=c_t|x|^{-(3-2t)}
\quad\text{for a.e. }x\neq0,
\]
which is \eqref{eq:Rieszbound}. Finally, since
\[
\frac{2}{(2_t^*)'}+\frac{3-2t}{3}=2.
\]
the Hardy--Littlewood--Sobolev inequality yields \eqref{eq:HLS}, with a constant independent of $a>0$.
\end{proof}

Throughout the sequel, we use the notation
\[
q_t=\frac{6}{3+2t},
\qquad
r_t=2q_t=\frac{12}{3+2t}.
\]
For $f,g\in L^{q_t}(\R^3)$, define
\[
B_a(f,g)=\int_{\R^3}(K_{t,a}*f)g\,dx.
\]
The evenness of $K_{t,a}$ implies that $B_a$
is symmetric. Moreover, \eqref{eq:HLS} yields
\begin{equation}\label{eq:Babound}
|B_a(f,g)|\leq C_t\|f\|_{q_t}\|g\|_{q_t},
\qquad
\|K_{t,a}*f\|_{q_t'}
\leq C_t\|f\|_{q_t}.
\end{equation}
Since $0<t<1$ and $4s+2t>3$, we have
$2<r_t<2_s^*$. Therefore, Lemma~\ref{lem:2.1}
and the equivalence of the $H^s$ and $H^s$
norms give
\begin{equation}\label{eq:u2est}
\|u\|_{r_t}\leq C\|u\|_{H^s},
\qquad
\|u^2\|_{q_t}
=\|u\|_{r_t}^2
\leq C\|u\|_{H^s}^2.
\end{equation}

\subsection{The electrostatic potential}\label{subsec:potential}

The estimate \eqref{eq:u2est}, together with the
continuous embedding of $\Dta$ into $L^{q_t'}(\R^3)$,
ensures that the source term $u^2$ defines a bounded
linear functional on $\Dta$. Thus the electrostatic
equation admits a unique potential for each
$u\in\Hs$, as shown below.

\begin{lemma}\label{lem:2.3}
For every $u \in \Hs$ there exists a unique $\varphi \in \Dta$ such that
\begin{equation}\label{eq:weakpot}
(\varphi,\psi)_{\Dta} = \int_{\R^3} u^2\psi\,dx \quad \text{for every } \psi \in \Dta,
\end{equation}
and it is given by
\[
\varphi = K_{t,a}*u^2 \quad \text{in } L^{q_t'}(\R^3).
\]
Moreover
\begin{equation}\label{eq:Taestimate}
\int_{\R^3} \varphi u^2\,dx = \|\varphi\|^2_{\Dta} \le C\|u\|^4_{H^s}.
\end{equation}
\end{lemma}

\begin{proof}
By \eqref{eq:u2est}, $u^2\in L^{q_t}(\R^3)$.
Moreover, since $\Dta\hookrightarrow L^{q_t'}(\R^3)$,
H\"older's inequality gives
\[
\left|\int_{\R^3}u^2\psi\,dx\right|
\leq \|u^2\|_{q_t}\|\psi\|_{q_t'}
\leq C\|u\|_{H^s}^2\|\psi\|_{\Dta},
\qquad \psi\in\Dta.
\]
Thus the right-hand side of \eqref{eq:weakpot}
defines a bounded linear functional on $\Dta$.
The Lax--Milgram theorem yields a unique
$\varphi_u\in\Dta$ satisfying \eqref{eq:weakpot},
with
\[
\|\varphi_u\|_{\Dta}
\leq C\|u\|_{H^s}^2.
\]
Taking $\psi=\varphi_u$ in \eqref{eq:weakpot}
gives
\[
\|\varphi_u\|_{\Dta}^2
=\int_{\R^3}\varphi_u u^2\,dx,
\]
and hence \eqref{eq:Taestimate}.

For the representation, choose
$f_k\in C_c^\infty(\R^3)$ such that
$f_k\to u^2$ in $L^{q_t}(\R^3)$, and set
$\varphi_k=K_{t,a}*f_k$. By the Fourier multiplier
formula,
\[
\widehat{\varphi_k}(\xi)
=\frac{\widehat f_k(\xi)}
{|\xi|^{2t}+a^{2t}|\xi|^{4t}},
\]
and
\[
\|\varphi_k\|_{\Dta}^2
=\int_{\R^3}
\frac{|\widehat f_k(\xi)|^2}
{|\xi|^{2t}+a^{2t}|\xi|^{4t}}\,d\xi
<\infty.
\]
The integral is finite since $2t<3$ and $\widehat f_k$ is rapidly
decreasing. Smooth annular cutoffs of $\widehat\varphi_k$ converge in
the displayed weighted norm by dominated convergence. Their inverse
transforms are Schwartz functions. Cutting off a Schwartz function in
physical space approximates it in $H^2$, hence in $\Dta$ since $2t<2$.
This places $\varphi_k$ in the specified completion and gives
\[
(\varphi_k,\psi)_{\Dta}
=\int_{\R^3}f_k\psi\,dx,
\qquad \psi\in\Dta.
\]
Subtracting \eqref{eq:weakpot} and taking
$\psi=\varphi_k-\varphi_u$, we obtain
\[
\begin{aligned}
\|\varphi_k-\varphi_u\|_{\Dta}^2
&=\int_{\R^3}(f_k-u^2)(\varphi_k-\varphi_u)\,dx\\
&\leq C\|f_k-u^2\|_{q_t}
\|\varphi_k-\varphi_u\|_{\Dta}.
\end{aligned}
\]
Consequently,
\[
\varphi_k\to\varphi_u
\quad\text{in }\Dta.
\]
On the other hand, by \eqref{eq:Babound},
\[
\|K_{t,a}*f_k-K_{t,a}*u^2\|_{q_t'}
\leq C_t\|f_k-u^2\|_{q_t}\to0.
\]
Since $\Dta\hookrightarrow L^{q_t'}(\R^3)$,
we also have $\varphi_k\to\varphi_u$ in
$L^{q_t'}(\R^3)$. Therefore,
\[
\varphi_u=K_{t,a}*u^2
\quad\text{in }L^{q_t'}(\R^3).
\]\end{proof}

The unique potential $\varphi_u$ allows us to
define the electrostatic energy associated with $u$ by
\[
T_a(u):=\int_{\R^3}\varphi_u u^2\,dx
=\|\varphi_u\|_{\Dta}^2
=B_a(u^2,u^2).
\]
By \eqref{eq:Taestimate},
\[
0\leq T_a(u)\leq C\|u\|_{H^s}^4.
\]
We also use $K_{t,0}=R_t$, $T_0(u)=\int(R_t*u^2)u^2$ and the
functional $J_{0,q}$ obtained by replacing $T_a$ by $T_0$ in
\eqref{eq:Ja}. Lemma~\ref{lem:2.2} and \eqref{eq:u2est} give
\begin{equation}\label{eq:uniformT}
0\le T_a(u)\le T_0(u)\le C\|u\|_{r_t}^4
\le C'\|u\|_{H^s}^4\qquad(a>0),
\end{equation}
where both constants are independent of $a$.
We next record the continuity of the electrostatic
potential with respect to the matter component.

\begin{lemma}\label{lem:2.4}
If $u_n^2 \to u^2$ in $L^{q_t}(\R^3)$, then
\[
\varphi_{u_n} \to \varphi_u \text{ in } \Dta, \qquad T_a(u_n) \to T_a(u).
\]
In particular this holds whenever $u_n \to u$ in $\Hs$.
\end{lemma}

\begin{proof}
Subtracting the equations \eqref{eq:weakpot} for $\varphi_{u_n}$ and $\varphi_u$ and testing with $\varphi_{u_n}-\varphi_u$ gives
\[
\|\varphi_{u_n}-\varphi_u\|^2_{\Dta} = \int_{\R^3} (u_n^2-u^2)(\varphi_{u_n}-\varphi_u)\,dx
\le \|u_n^2-u^2\|_{q_t}\|\varphi_{u_n}-\varphi_u\|_{q_t'} \to 0,
\]
since $\{\varphi_{u_n}\}$ is bounded in $\Dta \hookrightarrow L^{q_t'}$ by \eqref{eq:Taestimate}. The convergence of $T_a$ follows from $T_a(u) = \|\varphi_u\|^2_{\Dta}$.

If $u_n \to u$ in $\Hs$, then $u_n \to u$ in $L^{r_t}$ and
\[
\|u_n^2 - u^2\|_{q_t} \le \bigl( \|u_n\|_{r_t} + \|u\|_{r_t} \bigr)\|u_n-u\|_{r_t} \to 0,
\]
so the hypothesis is satisfied.
\end{proof}

\subsection{The variational functional}
\label{subsec:functional}

Substituting the unique potential from Lemma~\ref{lem:2.3} into the
matter equation gives the variational functional
\begin{equation}\label{eq:Ja}
J_{a,q}(u)
=\frac12\|u\|_{H^s}^2
+\frac{q^2}{4}T_a(u)
-\frac1p\int_{\R^3}|u|^p\,dx,
\qquad u\in\Hs.
\end{equation}
For fixed $a,q>0$, we simply write $J$ for
$J_{a,q}$ whenever no confusion can arise.

\begin{lemma}\label{lem:2.6}
The functional $J$ belongs to $C^1(\Hs,\R)$, and
\[
T_a'(u)[z]
=4B_a(u^2,uz)
=4\int_{\R^3}\varphi_uuz\,dx.
\]
Moreover,
\begin{equation}\label{eq:Jder}
\langle J'(u),v\rangle
=\int_{\R^3}
\left(
(-\Delta)^{s/2}u\,(-\Delta)^{s/2}v
+Vuv+q^2\varphi_uuv
-|u|^{p-2}uv
\right)\,dx.
\end{equation}
In particular, $J'(u)=0$ if and only if
$(u,\varphi_u)$ is a weak solution of \eqref{eq:sys},
and every critical point satisfies
\begin{equation}\label{eq:Nehari}
\|u\|_{H^s}^2+q^2T_a(u)=\|u\|_p^p.
\end{equation}
\end{lemma}

\begin{proof}
Since $B_a$ is continuous and symmetric, and
$u\mapsto u^2$ is continuously differentiable
from $\Hs$ into $L^{q_t}(\R^3)$, we have
$T_a\in C^1(\Hs,\R)$, with
\[
T_a'(u)[z]=4B_a(u^2,uz).
\]
Differentiating \eqref{eq:Ja} gives
\eqref{eq:Jder}. The equivalence between critical
points and weak solutions follows from Lemma~\ref{lem:2.3},
while \eqref{eq:Nehari} follows by taking $v=u$
in \eqref{eq:Jder}.
\end{proof}

\begin{lemma}\label{lem:2.7}
If $u\in\Hsr$ is a critical point of $J$
restricted to $\Hsr$, then $J'(u)=0$
in $(\Hs)^*$.
\end{lemma}

\begin{proof}
For $R\in SO(3)$, set
\[
(R\cdot v)(x)=v(R^{-1}x).
\]
The action is isometric on $\Hs$ and, by the
rotation invariance of the equation and of $T_a$,
\[
J(R\cdot v)=J(v).
\]
Since $u\in\Hsr$, we have $R\cdot u=u$.
Differentiating the invariance at $u$ gives
\[
\langle J'(u),R\cdot v\rangle
=\langle J'(u),v\rangle.
\]
Averaging over $SO(3)$, we obtain
\[
\langle J'(u),v\rangle
=
\left\langle J'(u),
\int_{SO(3)}R\cdot v\,dR\right\rangle.
\]
The averaged function belongs to $\Hsr$.
Hence the right-hand side vanishes by the
assumed radial criticality. Therefore
$J'(u)=0$ in $(\Hs)^*$.
\end{proof}

\section{Existence of solutions in \texorpdfstring{$H^s(\mathbb R^3)$}{Hs(R3)}}
\label{sec:existence}
In this section, we prove Theorems~\ref{thm:1.1}--\ref{thm:1.3}
in $H^s(\mathbb R^3)$. The first two results yield bounded
Palais--Smale sequences by different variational arguments,
whereas the third uses a truncated functional. In all three
cases, the nontriviality of a solution is obtained through
the same splitting lemma, whose proof is given in Appendix~A.

\subsection{A splitting lemma}
\begin{lemma}\label{lem:splitting}
Let $\{u_n\}\subset\Hs$ be a bounded Palais--Smale
sequence for $J_{a,q}$ at a level $d>0$, with
$u_n\rightharpoonup u^0$ in $\Hs$. Then $J_{a,q}'(u^0)=0$.
Up to a subsequence, either $u_n\to u^0$ strongly in $\Hs$, or
there exist $\ell\in\mathbb N$, nonzero critical points
$w^1,\ldots,w^\ell\in\Hs$ of $J_{a,q}$, and sequences
$\{y_n^j\}\subset\mathbb R^3$, $1\le j\le\ell$, such that
\begin{enumerate}[label=(\roman*)]
\item $|y_n^j|\to\infty$ for all $1\le j\le\ell$, and
$|y_n^i-y_n^j|\to\infty$ for all $i\ne j$;

\item
\[
\left\|u_n-u^0-\sum_{j=1}^{\ell}
w^j(\cdot-y_n^j)\right\|_{H^s}\to0;
\]

\item
\[
d=J_{a,q}(u^0)+\sum_{j=1}^{\ell}J_{a,q}(w^j);
\]

\item
\[
\|u_n\|_{H^s}^2=\|u^0\|_{H^s}^2
+\sum_{j=1}^{\ell}\|w^j\|_{H^s}^2+o(1).
\]
\end{enumerate}
\end{lemma}

\subsection{The case \texorpdfstring{$4\leq p<2_s^*$}{4 <= p < 2s*}}

\begin{proof}[Proof of Theorem~\ref{thm:1.1}]
Assume that $s>3/4$ and $4\leq p<2_s^*$.
Since $T_a(u)\geq0$ and $H^s(\mathbb R^3)\hookrightarrow
L^p(\mathbb R^3)$, there exist $\rho,\alpha>0$ such that
\[
J_{a,q}(u)\geq\alpha
\quad\text{whenever }\|u\|_{H^s}=\rho.
\]
To find a point with negative energy, choose
$w\in C_c^\infty(\mathbb R^3)\setminus\{0\}$ and set
$w_\tau(x)=\tau w(\tau x)$ for $\tau\geq1$. We have
\[
\|(-\Delta)^{s/2}w_\tau\|_2^2
=\tau^{2s-1}\|(-\Delta)^{s/2}w\|_2^2,\quad
\|w_\tau\|_2^2=\tau^{-1}\|w\|_2^2,\quad
\|w_\tau\|_p^p=\tau^{p-3}\|w\|_p^p.
\]
Moreover, the kernel estimate
$K_{t,a}(x)\leq C|x|^{-(3-2t)}$ gives
\[
T_a(w_\tau)\leq C\tau^{1-2t}.
\]
Since
\[
p-3\geq1>2s-1,\qquad p-3>-1,\qquad
p-3>1-2t,
\]
 we have $J_{0,q}(w_\tau)\to-\infty$ as well, because the same
Riesz bound holds at $a=0$. Choose $e=w_{\tau_0}$ with
$\|e\|_{H^s}>\rho$ and $J_{0,q}(e)<0$. This is possible since
$2s-1>0$. By \eqref{eq:uniformT}, $J_{a,q}(e)\le J_{0,q}(e)$ for
all $a>0$. Accordingly, we define
\[
\Gamma_e=\{\eta\in C([0,1],\Hs):\eta(0)=0,\ \eta(1)=e\},
\qquad
c_{a,q}=\inf_{\eta\in\Gamma_e}
\max_{0\leq\sigma\leq1}J_{a,q}(\eta(\sigma)).
\]
Then $c_{a,q}\geq\alpha>0$, and the mountain-pass theorem
provides a sequence $\{u_n\}\subset H^s(\mathbb R^3)$ such that
\[
J_{a,q}(u_n)\to c_{a,q},\qquad
J_{a,q}'(u_n)\to0\quad\text{in }(H^s)^*.
\]

We first observe that this sequence is bounded. Indeed,
since $p\geq4$ and $T_a(u_n)\geq0$,
\[
\begin{aligned}
pJ_{a,q}(u_n)-\langle J_{a,q}'(u_n),u_n\rangle
&=\left(\frac p2-1\right)\|u_n\|_{H^s}^2
+\left(\frac p4-1\right)q^2T_a(u_n)\\
&\geq\left(\frac p2-1\right)\|u_n\|_{H^s}^2.
\end{aligned}
\]
The left-hand side is
$O(1)+o(1)\|u_n\|_{H^s}$, so $\{u_n\}$ is bounded
in $H^s(\mathbb R^3)$.

After passing to a subsequence, let
$u_n\rightharpoonup u^0$ in $H^s(\mathbb R^3)$.
By Lemma~\ref{lem:splitting}, either $u_n\to u^0$
strongly in $H^s(\mathbb R^3)$, or there exist nonzero
critical points of $J_{a,q}$. In the first case,
\[
J_{a,q}(u^0)=\lim_{n\to\infty}J_{a,q}(u_n)
=c_{a,q}>0,
\]
and hence $u^0\ne0$. In the second case, any one of
the nonzero critical points gives a nontrivial solution.
This proves Theorem~\ref{thm:1.1}.
\end{proof}

\subsection{The case
\texorpdfstring{$p_*(s,t)<p<\min\{4,2_s^*\}$}
{$p_*(s,t)<p<\min\{4,2_s^*\}$}}
\label{subsec:subquartic-full}

Assume $p_*(s,t)<p<\min\{4,2_s^*\}$. Set
\begin{equation}\label{eq:gamma}
\gamma=s+t,\qquad d=4s+2t-3>0,\qquad b=(s+t)p-3>d .
\end{equation}
In fact, with $r_t=12/(3+2t)$ one has
\[
p_*(s,t)-r_t=\frac{2td}{(s+t)(3+2t)}>0,
\qquad
2_s^*-p_*(s,t)=\frac{2sd}{(3-2s)(s+t)}>0 .
\]
Since $p_*(s,t)<4$, the stated subquartic interval is nonempty for
every admissible pair $(s,t)$.

We recall the following bounded Palais--Smale result
of Guida and Rolando
 \cite[Theorem~1.1 and Remark~2.2]{GuidaRolando}.
\begin{lemma}\label{lem:GR-full}
Let $X$ be a Banach space and let
$I\in C^1(X,\mathbb R)$ be of the form
\[
I(u)=A(u)-B(u),
\]
where $A,B\in C^1(X,[0,\infty))$ and
$A(0)=B(0)=0$.
Assume that the following conditions are satisfied:

\begin{enumerate}[label={\rm(\roman*)}]
\item For every $0<\ell_1<\ell_2$, there exists
$\delta>0$ such that
\[
\operatorname{dist}(u,A_{\ell_1})\leq\delta
\quad\Longrightarrow\quad
u\in A_{\ell_2},
\qquad
A_\ell:=\{u\in X:A(u)\leq\ell\}.
\tag{A}
\]

\item There exists a sequence
$\{\psi_n\}\subset C(X,X)$ and two sequences
$\{\alpha_n\}$ and $\{\beta_n\}$ with
$\alpha_n>\beta_n>0$ such that, for all $u\in X$,
\[
A(u)\geq\alpha_n A(\psi_n(u)),
\qquad
B(u)\leq\beta_n B(\psi_n(u)),
\tag{$\Psi_1$}
\]
and
\[
\alpha_n\to1,\qquad \beta_n\to1,\qquad
\liminf_{n\to\infty}
\frac{|1-\beta_n|}{\alpha_n-\beta_n}<\infty.
\tag{$\Psi_2$}
\]

\item There exist $r>0$ and $e\in X$, with $\|e\|>r$,
such that
\[
\inf_{\|u\|=r}I(u)>I(0)\geq I(e),
\]
and
\[
\|\psi_n(0)\|\to0,\qquad
\|\psi_n(e)-e\|\to0.
\]
\end{enumerate}

Then there exists a sequence $\{u_n\}\subset X$ such that
\[
I(u_n)\to c,\qquad
I'(u_n)\to0\quad\text{in }X^*,
\qquad
\sup_n A(u_n)<\infty,
\]
where
\[
c:=\inf_{\eta\in\Gamma_e}
\max_{t\in[0,1]}I(\eta(t)),
\qquad
\Gamma_e:=
\{\eta\in C([0,1],X):
\eta(0)=0,\ \eta(1)=e\}.
\]
 The sequence can be chosen with $A(u_n)\le M$ for every $n$,
for any fixed
\begin{equation}\label{eq:grbound}
M>|c|\liminf_{n\to\infty}\frac{|1-\beta_n|}{\alpha_n-\beta_n}.
\end{equation}
\end{lemma}
\begin{proof}[Proof of Theorem~\ref{thm:1.2}]
We apply Lemma~\ref{lem:GR-full} to
$X=H^s(\mathbb R^3)$ and
$I=J_{a,q}$.
Writing
\[
J_{a,q}(u)=A_{a,q}(u)-B(u),
\]
where
\[
A_{a,q}(u)
=\frac12\|u\|_{H^s}^2+\frac{q^2}{4}T_a(u),
\qquad
B(u)=\frac1p\|u\|_p^p,
\]
we verify conditions {\rm (i)}--{\rm (iii)}
of Lemma~\ref{lem:GR-full}.

We first verify the mountain-pass geometry.
Since $A_{a,q}(u)\ge\frac12\|u\|_{H^s}^2$ and the Sobolev embedding gives
$B(u)\le C\|u\|_{H^s}^p$, there exist $r,\alpha>0$ such that
$J(u)\ge\alpha>0$ whenever $\|u\|_{H^s}=r$.
For $u_\tau(x)=\tau^{s+t}u(\tau x)$, changing variables gives
\[
\|(-\Delta)^{s/2}u_\tau\|_2^2=\tau^d\|(-\Delta)^{s/2}u\|_2^2,
\quad \|u_\tau\|_2^2=\tau^{d-2s}\|u\|_2^2,
\quad \|u_\tau\|_p^p=\tau^b\|u\|_p^p .
\]
For $u\in C_c^\infty$,
$\widehat{u_\tau^2}(\xi)=\tau^{2(s+t)-3}\widehat{u^2}(\xi/\tau)$.
Plancherel's identity for the field gives, for $\tau\ge1$,
\begin{equation}\label{eq:exactscaling}
\begin{aligned}
T_a(u_\tau)
 &=\tau^d\int_{\R^3}
 \frac{|\widehat{u^2}(\xi)|^2}{|\xi|^{2t}(1+(a\tau|\xi|)^{2t})}\,d\xi\\
 &=\tau^d T_{a\tau}(u)\le\tau^d T_a(u).
\end{aligned}
\end{equation}
At $a=0$ there is equality $T_0(u_\tau)=\tau^d T_0(u)$.
For fixed $a,\tau$, density and the continuity of $T_a,T_{a\tau},T_0$
extend these identities and inequalities to $u\in\Hs$. In particular this
argument covers $2_s^*\le4$.
Fix a nonzero function $w\in C_c^\infty(\mathbb R^3)$.Since $b>d>0$,
$J_{0,q}(w_\tau)\to-\infty$ and $\|w_\tau\|_{H^s}\to\infty$.
Choose $e=w_{\tau_0}$ with $\|e\|_{H^s}>r$ and
$J_{0,q}(e)<0$. The same $e$ works for all $a>0$ by \eqref{eq:uniformT}.
Then
\[
\Gamma_e:=\{\eta\in C([0,1],\Hs):\eta(0)=0,\ \eta(1)=e\},
\qquad
c_{a,q}:=\inf_{\eta\in\Gamma_e}\max_{\tau\in[0,1]}J(\eta(\tau))
\ge\alpha>0.
\]

We next verify condition {\rm (i)}.
For $\ell>0$, let
\[
A_\ell:=\{u\in\Hs:A_{a,q}(u)\leq\ell\}.
\]
If $v\in A_{\ell_1}$, then
$\|v\|_{H^s}\leq\sqrt{2\ell_1}$. By
Lemma~\ref{lem:2.6} and the continuity estimate for the
bilinear form defining $T_a$, we have
\[
\|T_a'(w)\|_{(\Hs)^*}\leq C_a\|w\|_{H^s}^3,
\qquad
\|A_{a,q}'(w)\|_{(\Hs)^*}
\leq C_{a,q}(\|w\|_{H^s}+\|w\|_{H^s}^3).
\]
Fix $0<\ell_1<\ell_2$. If
$\operatorname{dist}(u,A_{\ell_1})\leq\delta\leq1$,
choose $v\in A_{\ell_1}$ such that
$\|u-v\|_{H^s}<2\delta$. The segment $[v,u]$ is contained
in the ball of radius $\sqrt{2\ell_1}+2$; hence
\[
A_{a,q}(u)\leq A_{a,q}(v)
+C_{\ell_1,a,q}\|u-v\|_{H^s}
\leq\ell_1+2C_{\ell_1,a,q}\delta.
\]
Taking $\delta>0$ sufficiently small so that
$2C_{\ell_1,a,q}\delta<\ell_2-\ell_1$ yields
$u\in A_{\ell_2}$. Thus condition {\rm (i)} holds.

We now consider the dilation maps $\psi_n$ in condition {\rm (ii)}.
 By \eqref{eq:gamma} and \eqref{eq:exactscaling}, the largest exponent
among the three positive terms is $\widehat A=d$, while
$\widehat B=b>d$. Choose $\tau_n\downarrow1$,$\tau_n>1$, and define
\[
\psi_n(u):=u_{\tau_n},\qquad
\alpha_n:=\tau_n^{-\widehat A},\qquad
\beta_n:=\tau_n^{-\widehat B}.
\]
The scaling estimates imply
\[
A_{a,q}(u_{\tau_n})\leq\tau_n^{\widehat A}A_{a,q}(u),
\qquad
B(u_{\tau_n})=\tau_n^{\widehat B}B(u),
\]
whence
\[
A_{a,q}(u)\geq\alpha_n A_{a,q}(\psi_n(u)),
\qquad
B(u)\leq\beta_n B(\psi_n(u)).
\]
Furthermore,
\[
\alpha_n>\beta_n>0,\qquad
\alpha_n,\beta_n\to1,\qquad
\frac{|1-\beta_n|}{\alpha_n-\beta_n}
\longrightarrow
\frac{\widehat B}{\widehat B-\widehat A}<\infty.
\]
Thus {\rm ($\Psi$)} holds. Moreover, for each $n$, the
dilation operator $u\mapsto u_{\tau_n}$ is continuous on
$H^s(\mathbb R^3)$, and hence
$\psi_n\in C(X,X)$.

It remains to check the convergence of $\psi_n$ at $0$ and $e$.
Since $\psi_n(0)=0$ and $\tau_n\to1$, the strong continuity
of dilations at $\tau=1$ in $H^s(\mathbb R^3)$ gives
\[
\|\psi_n(e)-e\|_{H^s}
=\|e_{\tau_n}-e\|_{H^s}\to0.
\]
Indeed, this follows first for $C_c^\infty(\mathbb R^3)$
by direct calculation and then for $H^s(\mathbb R^3)$
by density and the uniform boundedness of the dilation
operators for $\tau$ near $1$. Since the $H^s$ norm is
equivalent to the standard $H^s$ norm, the required convergence
properties hold.

Consequently, all the hypotheses of Lemma~\ref{lem:GR-full}
are satisfied. Hence there exists a bounded sequence
$\{u_n\}\subset H^s(\mathbb R^3)$ such that
\[
J_{a,q}(u_n)\to c_{a,q}\geq\alpha>0,\qquad
J_{a,q}'(u_n)\to0
\quad\text{in }(H^s(\mathbb R^3))^*.
\]
After passing to a subsequence, we may assume that
$u_n\rightharpoonup u^0$ in $H^s(\mathbb R^3)$.
By Lemma~\ref{lem:splitting}, either $u_n\to u^0$
strongly in $H^s(\mathbb R^3)$, or there exist
$\ell\geq1$ and nonzero critical points
$w^1,\ldots,w^\ell$ of $J_{a,q}$.
In the first case,
\[
J_{a,q}(u^0)=\lim_{n\to\infty}J_{a,q}(u_n)
=c_{a,q}>0,
\]
so $u^0\ne0$. In the second case, any $w^j$
is a nontrivial critical point. Thus $J_{a,q}$
has a nontrivial critical point, which proves
Theorem~\ref{thm:1.2}.
\end{proof}

\subsection{The case of small coupling}

We now prove Theorem~\ref{thm:1.3}. Under its assumptions,
we establish the existence of a nontrivial solution for
sufficiently small coupling parameter $q$, using a truncated
variational functional and the splitting lemma.

For $T>0$, we introduce the following truncated functional.
Choose a nonincreasing function
$\chi\in C^1([0,\infty),[0,1])$ satisfying
\[
\chi(\tau)=
\begin{cases}
1, & 0\leq\tau\leq1,\\
0, & \tau\geq2,
\end{cases}
\qquad \chi'(\tau)\leq0.
\]
For $u\in H^s(\mathbb R^3)$, set
\[
\chi_T(u)=\chi\left(\frac{\|u\|_{H^s}^2}{T^2}\right)
\]
and define
\[
J_{a,q,T}(u)
=\frac12\|u\|_{H^s}^2
+\frac{q^2}{4}\chi_T(u)T_a(u)
-\frac1p\|u\|_p^p.
\]
Then $J_{a,q,T}\in C^1(H^s(\mathbb R^3),\mathbb R)$, and
\begin{align}
\langle J_{a,q,T}'(u),z\rangle
={}&(u,z)_{H^s}
+\frac{q^2}{4}\chi_T(u)T_a'(u)[z]
\label{eq:Jcutder}\\
&+\frac{q^2}{2T^2}
\chi'\left(\frac{\|u\|_{H^s}^2}{T^2}\right)
T_a(u)(u,z)_{H^s}
-\int_{\mathbb R^3}|u|^{p-2}uz\,dx.
\nonumber
\end{align}

\begin{lemma}\label{lem:cut-geometry}
Assume $2<p<2_s^*$. For every $T>0$ and $q>0$,
there exist $\rho,\alpha>0$ with $\rho<T$ such that
\[
J_{a,q,T}(u)\geq\alpha
\quad\text{if }\|u\|_{H^s}=\rho.
\]
Moreover, for every fixed $w\ne0$,
\[
J_{a,q,T}(\lambda w)\to-\infty
\quad\text{as }\lambda\to\infty.
\]
\end{lemma}

\begin{proof}
Since $T_a(u)\geq0$ and $H^s\hookrightarrow L^p$,
there exists $C_p>0$ such that
\[
J_{a,q,T}(u)\geq
\frac12\|u\|_{H^s}^2
-\frac{C_p}{p}\|u\|_{H^s}^p.
\]
The right-hand side is positive on a sufficiently small
sphere, and we choose its radius $\rho<T$.

On the other hand, for every fixed $w\ne0$, if $\lambda$
is sufficiently large so that
\[
\lambda^2\|w\|_{H^s}^2\geq2T^2,
\]
then $\chi_T(\lambda w)=0$ and
\[
J_{a,q,T}(\lambda w)
=\frac{\lambda^2}{2}\|w\|_{H^s}^2
-\frac{\lambda^p}{p}\|w\|_p^p
\longrightarrow-\infty
\quad\text{as }\lambda\to\infty.
\]
This establishes the mountain-pass geometry.
\end{proof}

In view of Lemma~\ref{lem:cut-geometry}, we define
\[
\Gamma_{a,q,T}
=
\left\{
\gamma\in C([0,1],H^s(\mathbb R^3)):
\gamma(0)=0,\ J_{a,q,T}(\gamma(1))<0
\right\}
\]
and the mountain-pass level
\[
c_{a,q,T}
=
\inf_{\gamma\in\Gamma_{a,q,T}}
\max_{\tau\in[0,1]}J_{a,q,T}(\gamma(\tau)).
\]
By Lemma~\ref{lem:cut-geometry}, $\Gamma_{a,q,T}\ne\varnothing$ and
$c_{a,q,T}\geq\alpha>0$. The mountain-pass theorem then
yields a Palais--Smale sequence $\{u_n\}\subset H^s(\mathbb R^3)$
such that
\[
J_{a,q,T}(u_n)\to c_{a,q,T},
\qquad
J_{a,q,T}'(u_n)\to0
\quad\text{in }(H^s(\mathbb R^3))^*.
\]

We show next that, for a suitable choice of $T>0$ and
sufficiently small $q>0$, the Palais--Smale sequence
$\{u_n\}$ satisfies
\[
\limsup_{n\to\infty}\|u_n\|_{H^s}<T.
\]
Consequently, $\chi_T(u_n)=1$ for all sufficiently large
$n$, and hence $\{u_n\}$ is also a Palais--Smale sequence
for $J_{a,q}$ at level $c_{a,q,T}$.
\begin{lemma}\label{lem:cut-bound}
Fix $V,s,t,p$ as above. There exist $T>0$ sufficiently large and $q_0>0$
sufficiently small such that, for every $a>0$ and
every $q\in(0,q_0)$, every $(PS)_{c_{a,q,T}}$
sequence $\{u_n\}$ for $J_{a,q,T}$ satisfies
\[
\limsup_{n\to\infty}\|u_n\|_{H^s}<T.
\]
\end{lemma}

\begin{proof}
Put
\[
A(u)=\|u\|_{H^s}^2,\qquad
\mathcal T_a(u)=T_a(u),\qquad
c_p=\frac12-\frac1p>0.
\]
By \eqref{eq:uniformT}, $\mathcal T_a(u)\le CA(u)^2$ with $C$
independent of $a,T$ and $q$.
Since $\mathcal T_a'(u)[u]=4\mathcal T_a(u)$, we have
\begin{align}
J_{a,q,T}(u)-\frac1p\langle J_{a,q,T}'(u),u\rangle
={}&c_pA(u)
+q^2\left(\frac14-\frac1p\right)\chi_T(u)\mathcal T_a(u)
\nonumber\\
&-\frac{q^2}{2pT^2}
\chi'\left(\frac{A(u)}{T^2}\right)\mathcal T_a(u)A(u).
\end{align}
Set $d_p=\max\{\frac1p-\frac14,0\}$. Since $\chi'\leq0$
and $\chi_T(u)=0$ for $A(u)\geq2T^2$, it follows that
\[
\chi_T(u)\mathcal T_a(u)\leq2CT^2A(u),
\]
and hence
\begin{equation}\label{eq:cutps-estimate}
J_{a,q,T}(u)-\frac1p\langle J_{a,q,T}'(u),u\rangle
\geq\bigl(c_p-2Cd_pq^2T^2\bigr)A(u).
\end{equation}
Write $\delta=c_p-2Cd_pq^2T^2$.

Fix a nonzero function $w\in C_c^\infty(\mathbb R^3)$ and choose
$R>0$ such that
\[
I_0(Rw):=\frac{R^2}{2}\|w\|_{H^s}^2
-\frac{R^p}{p}\|w\|_p^p<0.
\]
Choose $q^*>0$ so that
$I_0(Rw)+(q^*)^2R^4T_0(w)/4<0$. For $T>R\|w\|_{H^s}$ and
$0<q\leq q^*$, the path $\eta(\sigma)=\sigma Rw$ lies where
$\chi_T=1$. Its endpoint has negative energy for every $a>0$, and
\[
\max_{0\leq\sigma\leq1}J_{a,q,T}(\eta(\sigma))
\leq C_0:=\frac{R^2}{2}\|w\|_{H^s}^2
+\frac{(q^*)^2R^4}{4}T_0(w).
\]
This gives an upper bound $C_0$ for the mountain-pass level
in $H^s(\mathbb R^3)$. Here $q^*$ and $C_0$ depend only on the fixed data
and on $w,R$.

Choose $T$ sufficiently large such that
\[
T>R\|w\|_{H^s},\qquad c_pT^2>C_0,
\]
and then choose $q_0\leq q^*$ sufficiently small so that
\[
2Cd_pq_0^2T^4<c_pT^2-C_0.
\]
(The latter condition is automatic if $d_p=0$.)
Thus, for $0<q<q_0$, we have
\[
\delta>0,\qquad \delta T^2>C_0.
\]

Let $\{u_n\}$ be a $(PS)_{c_{a,q,T}}$ sequence for
$J_{a,q,T}$. By \eqref{eq:cutps-estimate},
\[
\delta\|u_n\|_{H^s}^2
\leq c_{a,q,T}+o(1)+o(1)\|u_n\|_{H^s}.
\]
Since $\delta>0$, the sequence $\{u_n\}$ is bounded in
$H^s(\mathbb R^3)$. It follows that
\[
\limsup_{n\to\infty}\|u_n\|_{H^s}^2
\leq\frac{c_{a,q,T}}{\delta}
\leq\frac{C_0}{\delta}<T^2.
\]
Hence, for all sufficiently large $n$,
\[
\chi_T(u_n)=1,\qquad
\chi'\left(\frac{\|u_n\|_{H^s}^2}{T^2}\right)=0.
\]
Consequently,
\[
J_{a,q,T}(u_n)=J_{a,q}(u_n),\qquad
J_{a,q,T}'(u_n)=J_{a,q}'(u_n),
\]
and $\{u_n\}$ is also a $(PS)_{c_{a,q,T}}$ sequence for
$J_{a,q}$. This completes the proof.
\end{proof}

We emphasize that the constants $T$ and $q_0$ can be
chosen independently of $a>0$.
By Lemma~\ref{lem:cut-bound}, $J_{a,q}$ admits a bounded
$(PS)_{c_{a,q,T}}$ sequence in $H^s(\mathbb R^3)$.
We now apply the splitting lemma stated at the beginning of this section.

\begin{proof}[Proof of Theorem~\ref{thm:1.3}]
Let $\{u_n\}$ be the bounded Palais--Smale sequence
for $J_{a,q,T}$ at the level $c_{a,q,T}>0$ obtained
above. By Lemma~\ref{lem:cut-bound}, after discarding finitely
many terms, $\{u_n\}$ is also a Palais--Smale sequence
for $J_{a,q}$ at the same level.

Passing to a subsequence, we may assume that
$u_n\rightharpoonup u^0$ in $\Hs$.
By Lemma~\ref{lem:splitting}, either
$u_n\to u^0$ strongly in $\Hs$, or
$J_{a,q}$ possesses a nonzero critical point.
In the first case,
\[
J_{a,q}(u^0)
=\lim_{n\to\infty}J_{a,q}(u_n)
=c_{a,q,T}>0,
\]
and hence $u^0\ne0$. In the second case, the
nonzero critical point of $J_{a,q}$ gives a
nontrivial solution of system \eqref{eq:sys}.
Therefore, in either case, system \eqref{eq:sys}
admits a nontrivial solution. This completes
the proof of Theorem~\ref{thm:1.3}.
\end{proof}

\section{Radial solutions and their limit}\label{sec:radial}

The variational constructions used in Section~3 can
also be carried out in the radial space
$H^s_{\mathrm{r}}(\mathbb R^3)$. Indeed, the dilation
maps involved in these constructions preserve radiality,
and the mountain-pass paths can be chosen in the radial
space. Moreover, the corresponding boundedness estimates
for the Palais--Smale sequences remain valid. We therefore
focus on the compactness argument and present it in detail
for the case $p\geq4$.

Assume that $s>\frac34$ and $4\leq p<2_s^*$.
As in the proof of Theorem~\ref{thm:1.1}, the
mountain-pass construction can be carried out in
$H^s_{\mathrm{r}}(\mathbb R^3)$, yielding a bounded
Palais--Smale sequence $\{u_n\}$ for $J_{a,q}$ at
a positive level $c_{a,q}^{\mathrm{r}}>0$, with
\[
u_n\in H^s_{\mathrm{r}}(\mathbb R^3),\qquad
J_{a,q}(u_n)\to c_{a,q}^{\mathrm{r}},\qquad
J'_{a,q}(u_n)\to0
\quad\text{in }(H^s_{\mathrm{r}})^*.
\]
Passing to a subsequence, we may assume that
\[
u_n\rightharpoonup u
\quad\text{in }H^s_{\mathrm{r}}(\mathbb R^3).
\]
By the compact radial embeddings, we have
\[
u_n\to u\quad\text{in }L^p(\mathbb R^3)
\cap L^{r_t}(\mathbb R^3).
\]
In particular, $u_n^2\to u^2$ in $L^{q_t}(\mathbb R^3)$.
By Lemma~\ref{lem:2.4}, it follows that
\[
\varphi_{u_n}\to\varphi_u
\quad\text{in }\mathcal D_a^{t,2}(\mathbb R^3),
\]
and hence in $L^{q_t'}(\mathbb R^3)$.

Passing to the limit in the variational identity
$J'_{a,q}(u_n)[v]=o(1)$ for
$v\in H^s_{\mathrm{r}}(\mathbb R^3)$, we obtain
\[
J'_{a,q}(u)[v]=0
\quad\text{for all }v\in H^s_{\mathrm{r}}(\mathbb R^3).
\]
Therefore, by Lemma~\ref{lem:2.7},
\[
J'_{a,q}(u)=0
\quad\text{in }(H^s(\mathbb R^3))^*.
\]

It remains to prove that $u_n\to u$ strongly in
$H^s(\mathbb R^3)$. Taking $u_n-u$ as a test
function in $J'_{a,q}(u_n)-J'_{a,q}(u)$, we obtain
\[
\begin{aligned}
\|u_n-u\|_{H^s}^2
={}&\int_{\mathbb R^3}
\bigl(|u_n|^{p-2}u_n-|u|^{p-2}u\bigr)(u_n-u)\,dx\\
&-q^2\int_{\mathbb R^3}
\bigl(\varphi_{u_n}u_n-\varphi_u u\bigr)(u_n-u)\,dx
+o(1).
\end{aligned}
\]
The first integral tends to zero by the strong
convergence in $L^p(\mathbb R^3)$. For the second
integral, H\"older's inequality gives
\[
\begin{aligned}
&\left|\int_{\mathbb R^3}
(\varphi_{u_n}u_n-\varphi_u u)(u_n-u)\,dx\right|\\
&\quad\leq
\|\varphi_{u_n}-\varphi_u\|_{q_t'}
\|u_n\|_{r_t}\|u_n-u\|_{r_t}
+\|\varphi_u\|_{q_t'}
\|u_n-u\|_{r_t}^2
\longrightarrow0.
\end{aligned}
\]
Consequently,
\[
u_n\to u\quad\text{in }H^s(\mathbb R^3).
\]
By the continuity of $J_{a,q}$, we have
\[
J_{a,q}(u)=\lim_{n\to\infty}J_{a,q}(u_n)
=c_{a,q}^{\mathrm{r}}>0.
\]
Thus $u\neq0$, and $(u,\varphi_u)$ is a
nontrivial radial solution of system~\eqref{eq:sys}.

The same compactness argument applies to the
Palais--Smale sequences obtained in the other two
cases. Hence system~\eqref{eq:sys} also admits
nontrivial radial solutions under the assumptions
of Theorems~\ref{thm:1.2} and~\ref{thm:1.3}.

\subsection[The limit of the electrostatic regularization]{{The limit as \texorpdfstring{$a\to0^+$}{a tends to zero}}}
Two estimates enter the limit. The first concerns the linear field
equation with a varying source. The second selects a bounded family of
radial critical points  at fixed coupling.
\begin{lemma}\label{lem:4.1}
Let $a_n\to0^+$, let $(f_n)$ be bounded in $L^{6/(3+2t)}(\R^3)$, and let $f^0\in L^{6/(3+2t)}(\R^3)$.
Let $\varphi_n$ be the unique solution in
$\mathcal D_{a_n}^{t,2}(\mathbb R^3)$ of
$(-\Delta)^t\varphi_n+a_n^{2t}(-\Delta)^{2t}\varphi_n=f_n$,
and let $\varphi_0$ be the unique solution in
$\mathcal D^{t,2}(\mathbb R^3)$ of
$(-\Delta)^t\varphi_0=f^0$.

Then the following implications hold.
\begin{enumerate}[label=(\roman*)]
\item if $f_n\rightharpoonup f^0$ in $L^{6/(3+2t)}$, then $\varphi_n\rightharpoonup\varphi_0$ in $D^{t,2}(\R^3)$;
\item if $f_n\to f^0$ strongly in $L^{6/(3+2t)}$, then $\varphi_n\to\varphi_0$ strongly in $D^{t,2}(\R^3)$ and $a_n^t(-\Delta)^t\varphi_n\to0$ in $L^2(\R^3)$.
\end{enumerate}
\end{lemma}

\begin{proof}
Testing the equation for $\varphi_n$ with $\varphi_n$ gives
\begin{equation}\label{eq:lim-energy}
\|(-\Delta)^{t/2}\varphi_n\|_2^2+a_n^{2t}\|(-\Delta)^t\varphi_n\|_2^2
=\int_{\R^3}f_n\varphi_n\,dx.
\end{equation}
By H\"older's inequality and the fractional Sobolev inequality,
\[
\left|\int f_n\varphi_n\,dx\right|
\le\|f_n\|_{6/(3+2t)}\|\varphi_n\|_{6/(3-2t)}
\le C\|f_n\|_{6/(3+2t)}\|(-\Delta)^{t/2}\varphi_n\|_2.
\]
The boundedness of $(f_n)$ and \eqref{eq:lim-energy} imply
\begin{equation}\label{eq:lim-est}
\|(-\Delta)^{t/2}\varphi_n\|_2\le C,
\qquad a_n^t\|(-\Delta)^t\varphi_n\|_2\le C.
\end{equation}
Thus, up to a subsequence, $\varphi_n\rightharpoonup\varphi^*$ in $D^{t,2}(\R^3)$. For every $\psi\in C_c^\infty(\R^3)$,
\[
\int (-\Delta)^{t/2}\varphi_n(-\Delta)^{t/2}\psi\,dx
+a_n^{2t}\int (-\Delta)^t\varphi_n(-\Delta)^t\psi\,dx
=\int f_n\psi\,dx.
\]
The higher-order term tends to zero, since
\[
\left|a_n^{2t}\int(-\Delta)^t\varphi_n(-\Delta)^t\psi\,dx\right|
\le a_n^t\|a_n^t(-\Delta)^t\varphi_n\|_2\|(-\Delta)^t\psi\|_2\to0.
\]
If $f_n\rightharpoonup f^0$, passage to the limit gives
\[
\int(-\Delta)^{t/2}\varphi^*(-\Delta)^{t/2}\psi\,dx=\int f^0\psi\,dx.
\]
By density and uniqueness in $D^{t,2}(\R^3)$, $\varphi^*=\varphi_0$. Every weakly convergent subsequence has the same limit, proving (i).

Suppose now that $f_n\to f^0$ strongly. By (i), $\varphi_n\rightharpoonup\varphi_0$ in $D^{t,2}(\R^3)$ and hence in $L^{6/(3-2t)}$. Moreover,
\[
\int f_n\varphi_n\,dx-\int f^0\varphi_0\,dx
=\int(f_n-f^0)\varphi_n\,dx+\int f^0(\varphi_n-\varphi_0)\,dx\to0.
\]
Testing the limiting equation with $\varphi_0$ yields
\[
\int f^0\varphi_0\,dx=\|(-\Delta)^{t/2}\varphi_0\|_2^2.
\]
It follows from \eqref{eq:lim-energy} that
\[
\|(-\Delta)^{t/2}\varphi_n\|_2^2+a_n^{2t}\|(-\Delta)^t\varphi_n\|_2^2
\longrightarrow\|(-\Delta)^{t/2}\varphi_0\|_2^2.
\]
Weak lower semicontinuity and nonnegativity of the second term imply convergence of the first norms. Weak convergence together with convergence of norms gives $\varphi_n\to\varphi_0$ strongly in $D^{t,2}(\R^3)$. Subtracting the first term from the energy identity also gives $a_n^t(-\Delta)^t\varphi_n\to0$ in $L^2$. This proves (ii).
\end{proof}

\begin{lemma}\label{lem:4.2}
Fix $s,t,p,V,q$ as in one of Theorems~\ref{thm:1.1}--\ref{thm:1.3}.
For every $a>0$, there exists a nontrivial radial solution
$u_a$ of \eqref{eq:sys} such that
\[
\sup_{a>0}\|u_a\|_{H^s}<\infty.
\]
\end{lemma}

\begin{proof}
In the two arbitrary-coupling cases, the mountain-pass
constructions of Section~\ref{sec:existence} can also be
carried out in the radial space by choosing the function
used to obtain the negative endpoint to be radial. Let
\[
\Gamma_e^{\rm r}
=\{\eta\in C([0,1],\Hsr):\eta(0)=0,\ \eta(1)=e\}.
\]
The positive sphere estimate is independent of $a$, and
the path $\eta(\sigma)=\sigma e$ gives
\begin{equation}\label{eq:radialupper}
0<\alpha\le c_{a,q}^{\rm r}
:=\inf_{\eta\in\Gamma_e^{\rm r}}
\max_{0\le\sigma\le1}J_{a,q}(\eta(\sigma))
\le C_*:=\max_{0\le\sigma\le1}J_{0,q}(\sigma e).
\end{equation}
All the levels in this comparison are obtained using radial paths.
When $p\ge4$, radial compactness gives a critical point at this level,
and the Nehari identity yields
\[
C_*\ge c_{a,q}^{\rm r}
=\frac{p-2}{2p}\|u_a\|_{H^s}^2
 +\frac{p-4}{4p}q^2T_a(u_a)
\ge\frac{p-2}{2p}\|u_a\|_{H^s}^2 .
\]
This includes $p=4$.

In the subquartic range of Theorem~\ref{thm:1.2},
$\alpha_n=\tau_n^{-d}$ and $\beta_n=\tau_n^{-b}$ give
\[
\lim_{n\to\infty}\frac{|1-\beta_n|}{\alpha_n-\beta_n}
=\frac b{b-d}.
\]
Use \eqref{eq:grbound} on the radial space with
$M=1+C_*b/(b-d)$. For every $a>0$ choose its Palais--Smale sequence
$\{u_{a,k}\}$ so that $A_{a,q}(u_{a,k})\le M$ for every $k$.
The fixed-$a$ radial compactness argument gives a strong limit along a
subsequence and thus $\|u_a\|_{H^s}^2\le2M$.

Finally, for $0<q<q_0$, we consider the radial minimax
level of the truncated functional. Choosing a nonzero radial
function in the construction of Lemma~\ref{lem:cut-bound},
we obtain a radial mountain-pass path and an upper bound
$C_0^{\rm r}$ for the radial minimax level, where
$C_0^{\rm r}$ is independent of $a$. Estimate
\eqref{eq:cutps-estimate}, tested with the radial
Palais--Smale sequence itself, gives
\[
\|u_a\|_{H^s}^2\leq C_0^{\rm r}/\delta<T^2
\]
after strong convergence. The constants $T$, $q_0$ and
$C_0^{\rm r}$ in this construction are independent of $a$.
\end{proof}

Since $2<r_t<2_s^*$ and $2<p<2_s^*$, the radial compact
embeddings and the corresponding dual embeddings yield
the following result.
\begin{lemma}\label{lem:4.3}
Let $a_n\to0^+$, and let
$u_n\in H^s_{\mathrm{r}}(\mathbb R^3)$
be radial weak solutions of system \eqref{eq:sys}
with parameter $a_n$. If
\[
\sup_n\|u_n\|_{H^s}<\infty,
\]
then, up to a subsequence, there exists
$u_0\in H^s_{\mathrm{r}}(\mathbb R^3)$ such that
\[
u_n\to u_0
\quad\text{in }H^s_{\mathrm{r}}(\mathbb R^3).
\]
\end{lemma}

\begin{proof}
By the boundedness of $\{u_n\}$, up to a subsequence,
\[
u_n\rightharpoonup u_0
\quad\text{in }H^s_{\mathrm{r}}(\mathbb R^3),
\qquad
u_n\to u_0
\quad\text{in }L^{r_t}(\mathbb R^3)\cap L^p(\mathbb R^3).
\]
Since $r_t=2q_t$, we have
$u_n^2\to u_0^2$ in $L^{q_t}(\mathbb R^3)$.
By Lemma~\ref{lem:4.1}(ii),
\[
\varphi_{a_n,u_n}\to\varphi_{0,u_0}
\quad\text{in }\mathcal D^{t,2}(\mathbb R^3)
\hookrightarrow L^{2_t^*}(\mathbb R^3).
\]
Since $1/r_t'=1/2_t^*+1/r_t$, H\"older's inequality
gives $\varphi_{a_n,u_n}u_n\to\varphi_{0,u_0}u_0$
in $L^{r_t'}(\mathbb R^3)$. Moreover,
$|u_n|^{p-2}u_n\to |u_0|^{p-2}u_0$
in $L^{p'}(\mathbb R^3)$. By the corresponding continuous
embeddings of the dual spaces, these convergences also
hold in $H^{-s}_{\mathrm{r}}(\mathbb R^3)$.

Passing to the limit in the weak formulation for radial
test functions, we obtain
\[
L_su_0+q^2\varphi_{0,u_0}u_0
=|u_0|^{p-2}u_0
\quad\text{in }H^{-s}_{\mathrm{r}}(\mathbb R^3),
\]
where $L_s=(-\Delta)^s+V$. Consequently,
\[
\begin{aligned}
L_s(u_n-u_0)
={}&|u_n|^{p-2}u_n-|u_0|^{p-2}u_0\\
&-q^2\big(
\varphi_{a_n,u_n}u_n-\varphi_{0,u_0}u_0
\big)
\longrightarrow0
\end{aligned}
\]
in $H^{-s}_{\mathrm{r}}(\mathbb R^3)$.
Since $L_s:H^s_{\mathrm{r}}(\mathbb R^3)
\to H^{-s}_{\mathrm{r}}(\mathbb R^3)$ is an
isomorphism, it follows that
$u_n\to u_0$ in $H^s_{\mathrm{r}}(\mathbb R^3)$.
\end{proof}

\begin{proof}[Proof of Theorem~\ref{thm:1.4}]
For the solutions selected by Lemma~\ref{lem:4.2}, the
required uniform bound holds. More generally, let
$\{(u_n,\varphi_n)\}$ be any sequence of nontrivial radial
solutions with parameters $a_n\to0^+$ satisfying
\[
\sup_n\|u_n\|_{H^s}<\infty.
\]
By Lemma~\ref{lem:4.3}, there exist
$u_0\in H^s_r(\mathbb R^3)$ and a subsequence such that

\[
u_n\to u_0
\quad\text{in }H^s_{\mathrm{r}}(\mathbb R^3).
\]
Since $u_n\to u_0$ in $L^{r_t}(\mathbb R^3)$ and
$r_t=2q_t$, we have $u_n^2\to u_0^2$ in
$L^{q_t}(\mathbb R^3)$. By Lemma~\ref{lem:4.1}(ii),
\[
\varphi_n=\varphi_{a_n,u_n}\to\varphi_{0,u_0}
\quad\text{in }\mathcal D^{t,2}(\mathbb R^3)
\hookrightarrow L^{2_t^*}(\mathbb R^3).
\]The same lemma yields
$a_n^t\|(-\Delta)^t\varphi_n\|_2\to0$.
We now show that $(u_0,\varphi_{0,u_0})$ solves
the limiting system \eqref{eq:sp}. The strong
convergence of $u_n$ and $\varphi_n$ yields
\[
|u_n|^{p-2}u_n
\to |u_0|^{p-2}u_0
\quad\text{in }L^{p'}(\mathbb R^3),
\]
and, since $1/r_t'=1/2_t^*+1/r_t$,
\[
\varphi_nu_n\to\varphi_{0,u_0}u_0
\quad\text{in }L^{r_t'}(\mathbb R^3).
\]
Let $\psi\in C_c^\infty(\mathbb R^3)$. The weak
formulation of the first equation in \eqref{eq:sys}
is
\[
\begin{aligned}
&\int_{\mathbb R^3}
(-\Delta)^{s/2}u_n(-\Delta)^{s/2}\psi\,dx
+V\int_{\mathbb R^3}u_n\psi\,dx\\
&\qquad
+q^2\int_{\mathbb R^3}\varphi_nu_n\psi\,dx
=\int_{\mathbb R^3}|u_n|^{p-2}u_n\psi\,dx.
\end{aligned}
\]
Passing to the limit, we obtain
\[
\begin{aligned}
&\int_{\mathbb R^3}
(-\Delta)^{s/2}u_0(-\Delta)^{s/2}\psi\,dx
+V\int_{\mathbb R^3}u_0\psi\,dx\\
&\qquad
+q^2\int_{\mathbb R^3}\varphi_{0,u_0}u_0\psi\,dx
=\int_{\mathbb R^3}|u_0|^{p-2}u_0\psi\,dx.
\end{aligned}
\]
Moreover, by the definition of $\varphi_{0,u_0}$,
$(-\Delta)^t\varphi_{0,u_0}=u_0^2$.
Consequently, $(u_0,\varphi_{0,u_0})$ solves
the limiting system \eqref{eq:sp}.

It remains to verify that $u_0$ is nontrivial.
Testing the first equation of \eqref{eq:sys} with
$u_n$, we obtain
\[
\|u_n\|_{H^s}^2+q^2T_{a_n}(u_n)
=\|u_n\|_p^p.
\]
By the Sobolev embedding and the nonnegativity of
$T_{a_n}(u_n)$, it follows that
\[
\|u_n\|_{H^s}^2
\leq \|u_n\|_p^p
\leq C\|u_n\|_{H^s}^p,
\]
where $C>0$ is independent of $n$. Since each $u_n$
is nontrivial, we deduce that
$\|u_n\|_{H^s}\geq C^{-1/(p-2)}>0$.
The strong convergence in $H^s(\mathbb R^3)$
therefore implies that $u_0\ne0$.
This completes the proof.
\end{proof}

\appendix
\section{Proof of the splitting lemma}
\label{app:splitting}
Throughout this appendix, we work under the standing
assumptions $a,V,q>0$, $4s+2t>3$, and $2<p<2_s^*$.
For brevity, we write
\[
H:=H^s(\mathbb R^3),\qquad
\|u\|_H:=\|u\|_{H^s}.
\]
We retain the notation $q_t$, $r_t$, $B_a$ and $T_a$
introduced in Section~\ref{sec:prelim}. In particular,
$2<r_t<2_s^*$, and
\[
H\hookrightarrow L^{r_t}(\mathbb R^3)
\cap L^p(\mathbb R^3)
\]
continuously and locally compactly. We shall also
use the estimates and derivative formula established
in Section~\ref{sec:prelim}.

All translations are understood in the usual sense.
Since $V$ is constant and the electrostatic term is
defined through convolution, the functional
$J_{a,q}$ is invariant under translations.

\begin{lemma}\label{lem:A.1}
Let $\{u_n\}\subset H$ be bounded and suppose
\[
u_n\rightharpoonup u\quad\text{in }H,
\qquad \|J_{a,q}'(u_n)\|_{H^*}\to0.
\]
Then $J_{a,q}'(u)=0$.
\end{lemma}

\begin{proof}
After passing to a subsequence, local compactness gives
\[
u_n\to u\quad\text{in }L^{r_t}_{\rm loc}(\mathbb R^3)\cap L^p_{\rm loc}(\mathbb R^3),
\qquad u_n(x)\to u(x)\quad\text{a.e.}
\]
Since $\{u_n\}$ is bounded in $L^{r_t}$,
\[
u_n^2\rightharpoonup u^2\quad\text{in }L^{q_t}(\mathbb R^3).
\]
For $\psi\in C_c^\infty(\mathbb R^3)$, the derivative is
\[
J_{a,q}'(u_n)[\psi]
=(u_n,\psi)_{H^s}+q^2B_a(u_n^2,u_n\psi)
-\int_{\mathbb R^3}|u_n|^{p-2}u_n\psi\,dx.
\]
The first term converges by weak convergence in $H$. Since $u_n\psi\to u\psi$ in $L^{q_t}$, the HLS estimate and weak convergence of the squares imply
\[
B_a(u_n^2,u_n\psi)\to B_a(u^2,u\psi).
\]
The nonlinear term converges by local strong convergence in $L^p$. Hence
\[
J_{a,q}'(u)[\psi]=0\qquad\forall\psi\in C_c^\infty(\mathbb R^3).
\]
Density of $C_c^\infty(\mathbb R^3)$ in $H$ yields $J_{a,q}'(u)=0$ in $H^*$.
\end{proof}

\begin{lemma}\label{lem:A.2}
Let $u\in H$ and let $v_n\rightharpoonup0$ be bounded in $H$. Then
\[
J_{a,q}(u+v_n)=J_{a,q}(u)+J_{a,q}(v_n)+o(1),
\]
and
\[
\|J_{a,q}'(u+v_n)-J_{a,q}'(u)-J_{a,q}'(v_n)\|_{H^*}\to0.
\]
\end{lemma}

\begin{proof}
After passing to a subsequence,
\[
v_n\to0\quad\text{in }L^{r_t}_{\rm loc}\cap L^p_{\rm loc},
\qquad v_n(x)\to0\quad\text{a.e.}
\]
We first note that the multiplication map
\[
M_u:H\to L^{q_t}(\mathbb R^3),\qquad M_u z=uz,
\]
is compact. Indeed, approximate $u$ in $L^{r_t}$ by functions $u_m\in C_c^\infty$. Hölder's inequality gives
\[
\|(u-u_m)z\|_{q_t}
\le \|u-u_m\|_{r_t}\|z\|_{r_t}
\le C\|u-u_m\|_{r_t}\|z\|_H.
\]
For fixed $m$, compactness follows from local compactness of $H\hookrightarrow L^{r_t}_{\rm loc}$. Thus $M_u$ is a norm limit of compact operators. In particular,
\begin{equation}\label{eq:A-uv}
uv_n\to0\quad\text{in }L^{q_t}.
\end{equation}
Also, boundedness in $L^{r_t}$ and almost-everywhere convergence imply
\begin{equation}\label{eq:A-vsquare}
v_n^2\rightharpoonup0\quad\text{in }L^{q_t}.
\end{equation}

The Hilbert norm splits, and the Br\'ezis--Lieb lemma \cite{BrezisLieb} gives
\[
\|u+v_n\|_H^2=\|u\|_H^2+\|v_n\|_H^2+o(1),
\qquad
\|u+v_n\|_p^p=\|u\|_p^p+\|v_n\|_p^p+o(1).
\]
Expanding $T_a(u+v_n)=B_a((u+v_n)^2,(u+v_n)^2)$,
every term containing $uv_n$ tends to zero by
\eqref{eq:Babound} and \eqref{eq:A-uv}. The remaining
mixed term tends to zero since
\[
B_a(u^2,v_n^2)
=\int_{\mathbb R^3}\varphi_u v_n^2\,dx
\longrightarrow0,
\qquad
\varphi_u=K_{t,a}*u^2\in L^{q_t'},
\]
by \eqref{eq:A-vsquare}. Thus
\[
T_a(u+v_n)=T_a(u)+T_a(v_n)+o(1),
\]
which proves the energy splitting.

For $F(x)=|x|^{p-2}x$, the mean value estimate, applied first for
$|x|\ge|y|$ and then with $x,y$ interchanged, gives
\[
|F(x+y)-F(x)-F(y)|
\le C_p\bigl(|x|^{p-2}|y|+|x||y|^{p-2}\bigr).
\]
On any ball the right-hand side with $(x,y)=(u,v_n)$ tends to zero in
$L^{p'}$ by local $L^p$ convergence and H\"older's inequality. Outside
$B_R$ its norm is bounded by
$C_p(\|u\|_{L^p(B_R^c)}^{p-2}M+\|u\|_{L^p(B_R^c)}M^{p-2})$,
where $M=\sup_n\|v_n\|_p$. Letting $R\to\infty$ proves
\begin{equation}\label{eq:A-power-der}
\big\|F(u+v_n)-F(u)-F(v_n)\big\|_{p'}\to0.
\end{equation}
Consequently, by H\"older's inequality and the
embedding $H\hookrightarrow L^p$,
\[
\sup_{\|z\|_H\le1}
\left|\int_{\mathbb R^3}
\big(F(u+v_n)-F(u)-F(v_n)\big)z\,dx\right|
\to0.
\]
It remains to estimate the electrostatic part.
For every $z\in H$, we have
\begin{align}
&\frac14\big[T_a'(u+v_n)-T_a'(u)-T_a'(v_n)\big][z]\notag\\
&\quad=B_a(u^2,v_nz)+B_a(v_n^2,uz)
+2B_a(uv_n,uz)+2B_a(uv_n,v_nz).
\label{eq:A-electro-exp}
\end{align}
We claim that
\[
\sup_{\|z\|_H\le1}
\left|
B_a(u^2,v_nz)+B_a(v_n^2,uz)
+2B_a(uv_n,uz)+2B_a(uv_n,v_nz)
\right|\to0.
\]
For the first term, using
$\varphi_u=K_{t,a}*u^2\in L^{q_t'}$ and
$\frac1{r_t'}=\frac1{r_t}+\frac1{q_t'}$,
we have
\[
|B_a(u^2,v_nz)|
\le \|\varphi_uv_n\|_{r_t'}\|z\|_{r_t}.
\]
Local strong convergence of $v_n$ in $L^{r_t}$ and
the small $L^{q_t'}$-tail of $\varphi_u$ imply
\[
\|\varphi_uv_n\|_{r_t'}\to0.
\]
Consequently,
\[
\sup_{\|z\|_H\le1}|B_a(u^2,v_nz)|\to0.
\]

For the second term, set
$F_n(g)=B_a(v_n^2,g)$ on $L^{q_t}$.
By \eqref{eq:Babound} and the boundedness of $v_n$
in $L^{r_t}$, the functionals $F_n$ are uniformly
bounded. Moreover, \eqref{eq:A-vsquare} implies
that $F_n(g)\to0$ for every fixed $g\in L^{q_t}$.
Since $M_u:H\to L^{q_t}$ is compact, the set
\[
\mathcal K_u:=\{uz:\|z\|_H\le1\}
\]
is relatively compact in $L^{q_t}$. Hence the
uniform boundedness of $F_n$ and their pointwise
convergence to zero yield
\[
\sup_{\|z\|_H\le1}|B_a(v_n^2,uz)|\to0.
\]

Finally, \eqref{eq:Babound} and \eqref{eq:A-uv}
give
\[
\begin{aligned}
&|B_a(uv_n,uz)|+|B_a(uv_n,v_nz)|\\
&\qquad\le C\|uv_n\|_{q_t}
\big(\|uz\|_{q_t}+\|v_nz\|_{q_t}\big)\\
&\qquad\le o(1)\|z\|_H.
\end{aligned}
\]
Combining these estimates with \eqref{eq:A-power-der}
and \eqref{eq:A-electro-exp} proves the derivative
splitting. Since the argument applies to a further
subsequence of every subsequence, the assertions
hold for the original sequence.
\end{proof}

We recall the fractional vanishing lemma
of Felmer, Quaas and Tan \cite[Lemma~2.2]{FelmerQuaasTan}.
\begin{lemma}\label{lem:A.3}
Let $\{v_n\}\subset H$ be bounded and suppose that
\[
\sup_{y\in\mathbb R^3}
\int_{B_1(y)}|v_n|^2\,dx\to0.
\]
Then
\[
 v_n\to0\quad\text{in }L^r(\mathbb R^3)\quad\text{for every }2<r<2_s^*.
\]
\end{lemma}

\begin{lemma}\label{lem:A.4}
Suppose $v_n\rightharpoonup0$ in $H$ and
$J_{a,q}'(v_n)\to0$ in $H^*$. After passage to a subsequence, either
\[
v_n\to0\quad\text{strongly in }H,
\]
or there exist $y_n\in\mathbb R^3$ and
$w\in H\setminus\{0\}$ satisfying
\[
|y_n|\to\infty,\qquad
v_n(\cdot+y_n)\rightharpoonup w\quad\text{in }H,
\qquad J_{a,q}'(w)=0.
\]
\end{lemma}

\begin{proof}
If
\[
\sup_{y\in\mathbb R^3}\int_{B_1(y)}|v_n|^2\,dx\to0,
\]
Lemma~\ref{lem:A.3} gives $v_n\to0$ in $L^p(\mathbb R^3)$.
Testing $J_{a,q}'(v_n)$ against $v_n$ yields
\[
\|v_n\|_H^2+q^2T_a(v_n)-\|v_n\|_p^p=o(1).
\]
Since $T_a(v_n)\ge0$, it follows that $v_n\to0$ strongly in $H$.

Otherwise, after passing to a subsequence, there exist
$\delta>0$ and $y_n\in\mathbb R^3$ such that
\[
\int_{B_1(y_n)}|v_n|^2\,dx\ge\delta.
\]
Necessarily $|y_n|\to\infty$, since bounded centers would
contradict local compactness and $v_n\rightharpoonup0$.
Define $\widetilde v_n=v_n(\cdot+y_n)$. Translation invariance
gives boundedness and
\[
\|J_{a,q}'(\widetilde v_n)\|_{H^*}
=\|J_{a,q}'(v_n)\|_{H^*}\to0.
\]
After passing to a subsequence, $\widetilde v_n\rightharpoonup w$
in $H$. By the local compactness of the embedding
$H^s(\mathbb R^3)\hookrightarrow L^2_{\mathrm{loc}}(\mathbb R^3)$,
we have $\widetilde v_n\to w$ in $L^2(B_1)$. Hence
\[
\int_{B_1}|w|^2\,dx\ge\delta.
\]
Thus $w\ne0$. Lemma~\ref{lem:A.1} gives $J_{a,q}'(w)=0$.
\end{proof}

\begin{lemma}\label{lem:A.5}
Let $v\in H$ and let $\{v_n\}\subset H$ be bounded.
\begin{itemize}
\item[(i)] If $\{y_n\}\subset\mathbb R^3$ satisfies
$|y_n|\to\infty$, then
\[
v(\cdot+y_n)\rightharpoonup0\quad\text{in }H.
\]

\item[(ii)] If $\{y_n\}\subset\mathbb R^3$ is bounded and
$v_n\rightharpoonup0$ in $H$, then, up to a subsequence,
\[
v_n(\cdot+y_n)\rightharpoonup0\quad\text{in }H.
\]
\end{itemize}
\end{lemma}

\begin{proof}
(i) For any $\psi\in H$, the Fourier representation of the
inner product gives
\[
((v\cdot+y_n),\psi)_{H^s}
=\int_{\mathbb R^3}e^{iy_n\cdot\xi}
(|\xi|^{2s}+V)\widehat v(\xi)
\overline{\widehat\psi(\xi)}\,d\xi.
\]
Since the integrand without the exponential factor belongs
to $L^1(\mathbb R^3)$ by the Cauchy--Schwarz inequality,
the Riemann--Lebesgue lemma yields
$((v\cdot+y_n),\psi)_{H^s}\to0$.

(ii) Passing to a subsequence, we may assume that $y_n\to y$.
For any $\psi\in H$, we have
\[
\begin{aligned}
(v_n(\cdot+y_n),\psi)_{H^s}
&=(v_n,\psi(\cdot-y_n))_{H^s}\\
&=(v_n,\psi(\cdot-y))_{H^s}+o(1)\to0,
\end{aligned}
\]
where $o(1)$ follows from the strong continuity of translations
in $H$ and the boundedness of $\{v_n\}$, while the convergence
to zero follows from $v_n\rightharpoonup0$ in $H$.
This proves (ii).
\end{proof}

\begin{proof}[Proof of Lemma~\ref{lem:splitting}]
By boundedness, after passing to a subsequence,
\[
u_n\rightharpoonup u^0\quad\text{in }H.
\]
Lemma~\ref{lem:A.1} gives $J_{a,q}'(u^0)=0$. Set
\[
v_n^0:=u_n-u^0.
\]
Then $v_n^0\rightharpoonup0$ in $H$. Lemma~\ref{lem:A.2} gives
\[
J(u_n)=J(u^0)+J(v_n^0)+o(1),
\qquad \|J'(v_n^0)\|_{H^*}\to0,
\]
where $J=J_{a,q}$. Thus $v_n^0$ satisfies the hypotheses of Lemma~\ref{lem:A.4}.

\medskip
\noindent\textit{First profile.}
If $v_n^0\to0$ strongly in $H$, the extraction terminates.
Otherwise, Lemma~\ref{lem:A.4} yields a sequence
$\{y_n^1\}\subset\mathbb R^3$ and a nonzero function
$w^1\in H$ such that
\[
|y_n^1|\to\infty,\qquad
v_n^0(\cdot+y_n^1)\rightharpoonup w^1
\quad\text{in }H,
\qquad J'(w^1)=0.
\]
Set
\[
v_n^1:=v_n^0-w^1(\cdot-y_n^1).
\]
Then
\[
v_n^1(\cdot+y_n^1)
=v_n^0(\cdot+y_n^1)-w^1\rightharpoonup0
\quad\text{in }H.
\]
On the other hand, Lemma~\ref{lem:A.5}(i) and
$|y_n^1|\to\infty$ imply
\[
w^1(\cdot-y_n^1)\rightharpoonup0
\quad\text{in }H.
\]
Consequently, $v_n^1\rightharpoonup0$ in the original frame.

Since translations are isometries in $H$ and $J$ is
translation invariant, the sequence
$\{v_n^0(\cdot+y_n^1)\}$ is also a bounded
Palais--Smale sequence for $J$. Moreover,
\[
v_n^0(\cdot+y_n^1)
=w^1+v_n^1(\cdot+y_n^1),
\qquad
v_n^1(\cdot+y_n^1)\rightharpoonup0
\quad\text{in }H.
\]
Applying Lemma~\ref{lem:A.2} and using $J'(w^1)=0$,
we obtain
\[
\|J'(v_n^1)\|_{H^*}\to0,\qquad
J(v_n^0)=J(w^1)+J(v_n^1)+o(1).
\]
Moreover, $\{v_n^1\}$ is bounded in $H$. Since
$J(v_n^0)$ converges, the energy splitting shows that
$J(v_n^1)$ also converges. Thus $\{v_n^1\}$ is a
bounded Palais--Smale sequence for $J$.

If $v_n^1\to0$ strongly in $H$, then, by the definition
of $v_n^1$,
\[
u_n=u^0+w^1(\cdot-y_n^1)+o_H(1).
\]
Furthermore, since $v_n^0\rightharpoonup0$ in $H$, we have
\[
\|u_n\|_H^2
=\|u^0\|_H^2+\|v_n^0\|_H^2+o(1).
\]
Moreover, the identity
\[
v_n^0(\cdot+y_n^1)
=w^1+v_n^1(\cdot+y_n^1)
\]
and the strong convergence $v_n^1\to0$ in $H$ yield
\[
\|v_n^0\|_H^2=\|w^1\|_H^2+o(1).
\]
Consequently,
\[
\|u_n\|_H^2
=\|u^0\|_H^2+\|w^1\|_H^2+o(1).
\]
Moreover, by the energy splittings obtained above and
the continuity of $J$, we have
\[
J(u_n)=J(u^0)+J(w^1)+J(v_n^1)+o(1)
      =J(u^0)+J(w^1)+o(1).
\]
Since $J(u_n)\to d$, we obtain
\[
d=J(u^0)+J(w^1).
\]
Thus the conclusion of Lemma~\ref{lem:splitting} holds
with $\ell=1$.

\medskip
\noindent\textit{Second profile.}
If $v_n^1\not\to0$ strongly in $H$, then, since
$\{v_n^1\}$ is a bounded Palais--Smale sequence for $J$,
Lemma~\ref{lem:A.4} yields a sequence
$\{y_n^2\}\subset\mathbb R^3$ and a nonzero function
$w^2\in H$ such that
\[
|y_n^2|\to\infty,\qquad
v_n^1(\cdot+y_n^2)\rightharpoonup w^2
\quad\text{in }H,\qquad J'(w^2)=0.
\]
We first show that
\[
|y_n^2-y_n^1|\to\infty.
\]
Suppose otherwise. Then, up to a subsequence,
$\{y_n^2-y_n^1\}$ is bounded. Since
$v_n^1(\cdot+y_n^1)\rightharpoonup0$ in $H$,
Lemma~\ref{lem:A.5}(ii), applied to
$v_n^1(\cdot+y_n^1)$ with the translation sequence
$y_n^2-y_n^1$, yields, up to a further subsequence,
\[
v_n^1(\cdot+y_n^2)\rightharpoonup0
\quad\text{in }H.
\]
This contradicts
$v_n^1(\cdot+y_n^2)\rightharpoonup w^2\ne0$.
Therefore, $|y_n^2-y_n^1|\to\infty$.

Set
\[
v_n^2:=v_n^1-w^2(\cdot-y_n^2).
\]
Then
\[
v_n^2(\cdot+y_n^2)
=v_n^1(\cdot+y_n^2)-w^2\rightharpoonup0
\quad\text{in }H.
\]
Moreover, since $|y_n^2-y_n^1|\to\infty$,
Lemma~\ref{lem:A.5}(i) implies
\[
w^2(\cdot+y_n^1-y_n^2)\rightharpoonup0
\quad\text{in }H.
\]
Consequently,
\[
v_n^2(\cdot+y_n^1)
=v_n^1(\cdot+y_n^1)
-w^2(\cdot+y_n^1-y_n^2)
\rightharpoonup0
\quad\text{in }H.
\]
On the other hand, since $|y_n^2|\to\infty$,
Lemma~\ref{lem:A.5}(i) gives
\[
w^2(\cdot-y_n^2)\rightharpoonup0
\quad\text{in }H.
\]
Together with $v_n^1\rightharpoonup0$ in $H$, this yields
\[
v_n^2\rightharpoonup0
\quad\text{in }H.
\]

To show that $\{v_n^2\}$ is again a bounded
Palais--Smale sequence, set
\[
z_n^2:=v_n^1(\cdot+y_n^2).
\]
Since translations are isometries in $H$ and $J$ is
translation invariant, $\{z_n^2\}$ is a bounded
Palais--Smale sequence for $J$. Furthermore,
\[
z_n^2\rightharpoonup w^2,\qquad
z_n^2-w^2=v_n^2(\cdot+y_n^2)\rightharpoonup0
\quad\text{in }H.
\]
Applying Lemma~\ref{lem:A.2} and using $J'(w^2)=0$,
we obtain
\[
\|J'(z_n^2-w^2)\|_{H^*}\to0
\]
and
\[
J(z_n^2)=J(w^2)+J(z_n^2-w^2)+o(1).
\]
By translation invariance, it follows that
\[
\|J'(v_n^2)\|_{H^*}\to0,\qquad
J(v_n^1)=J(w^2)+J(v_n^2)+o(1).
\]
Moreover, $\{v_n^2\}$ is bounded in $H$. Since
$J(v_n^1)$ converges, the energy splitting shows that
$J(v_n^2)$ also converges. Thus $\{v_n^2\}$ is a
bounded Palais--Smale sequence for $J$.

If $v_n^2\to0$ strongly in $H$, then, by the definitions
of $v_n^1$ and $v_n^2$,
\[
u_n=u^0+w^1(\cdot-y_n^1)
+w^2(\cdot-y_n^2)+o_H(1).
\]
Furthermore, since $v_n^0\rightharpoonup0$ in $H$,
\[
\|u_n\|_H^2
=\|u^0\|_H^2+\|v_n^0\|_H^2+o(1).
\]
The identity
\[
v_n^0(\cdot+y_n^1)
=w^1+v_n^1(\cdot+y_n^1)
\]
and the weak convergence
$v_n^1(\cdot+y_n^1)\rightharpoonup0$ in $H$ give
\[
\|v_n^0\|_H^2
=\|w^1\|_H^2+\|v_n^1\|_H^2+o(1).
\]
Similarly, the identity
\[
v_n^1(\cdot+y_n^2)
=w^2+v_n^2(\cdot+y_n^2)
\]
and the weak convergence
$v_n^2(\cdot+y_n^2)\rightharpoonup0$ in $H$ yield
\[
\|v_n^1\|_H^2
=\|w^2\|_H^2+\|v_n^2\|_H^2+o(1).
\]
Since $v_n^2\to0$ strongly in $H$, we conclude that
\[
\|u_n\|_H^2
=\|u^0\|_H^2+\|w^1\|_H^2
+\|w^2\|_H^2+o(1).
\]
Moreover, combining the energy splittings obtained above
and using $J(v_n^2)\to J(0)=0$, we have
\[
J(u_n)
=J(u^0)+J(w^1)+J(w^2)+J(v_n^2)+o(1)
=J(u^0)+J(w^1)+J(w^2)+o(1).
\]
Since $J(u_n)\to d$, it follows that
\[
d=J(u^0)+J(w^1)+J(w^2).
\]
Thus the conclusion of Lemma~\ref{lem:splitting} holds
with $\ell=2$. Otherwise, we continue
the extraction procedure as described below.

Suppose that, for some $k\geq2$, $\{v_n^k\}$ is a
bounded Palais--Smale sequence for $J$ such that
\[
v_n^k\rightharpoonup0\quad\text{in }H,\qquad
v_n^k(\cdot+y_n^i)\rightharpoonup0
\quad\text{in }H,\quad i=1,\ldots,k.
\]
If $v_n^k\to0$ strongly in $H$, the extraction process
terminates. Otherwise, Lemma~A.4 yields a sequence
$\{y_n^{k+1}\}\subset\mathbb R^3$ and a nonzero function
$w^{k+1}\in H$ such that
\[
|y_n^{k+1}|\to\infty,\qquad
v_n^k(\cdot+y_n^{k+1})\rightharpoonup w^{k+1}
\quad\text{in }H,\qquad J'(w^{k+1})=0.
\]
Moreover, the translation parameters satisfy
\[
|y_n^{k+1}-y_n^i|\to\infty,
\qquad i=1,\ldots,k.
\]
Indeed, if this were false for some $i$, then, up to a
subsequence, $\{y_n^{k+1}-y_n^i\}$ would be bounded.
By Lemma~A.5(ii), applied to $v_n^k(\cdot+y_n^i)$, we
would obtain
\[
v_n^k(\cdot+y_n^{k+1})\rightharpoonup0
\quad\text{in }H,
\]
which contradicts $w^{k+1}\ne0$.

Set
\[
v_n^{k+1}:=v_n^k-w^{k+1}(\cdot-y_n^{k+1}).
\]
Then, by Lemma~A.5(i) and the separation of the
translation parameters, we have
\[
v_n^{k+1}(\cdot+y_n^i)\rightharpoonup0
\quad\text{in }H,\quad i=1,\ldots,k+1,
\]
and also $v_n^{k+1}\rightharpoonup0$ in $H$.
Furthermore, applying Lemma~A.2 in the translated
coordinates, as in the preceding steps, we obtain
\[
\|J'(v_n^{k+1})\|_{H^*}\to0,
\]
\[
J(v_n^k)=J(w^{k+1})+J(v_n^{k+1})+o(1),
\]
and
\[
\|v_n^k\|_H^2
=\|w^{k+1}\|_H^2+\|v_n^{k+1}\|_H^2+o(1).
\]
In addition, $\{v_n^{k+1}\}$ is bounded in $H$.
Since $J(v_n^k)$ converges, the energy splitting implies
that $J(v_n^{k+1})$ also converges. Thus
$\{v_n^{k+1}\}$ is again a bounded Palais--Smale
sequence for $J$, and all the induction hypotheses are
satisfied for the next remainder.

It remains to prove that only finitely many profiles can
be extracted. For any nontrivial critical point $w$,
testing $J'(w)=0$ with $w$ gives
\[
\|w\|_{H^s}^2+q^2T_a(w)
=\|w\|_p^p
\leq C\|w\|_{H^s}^p,
\]
where $C>0$ is independent of $w$. Since $T_a(w)\geq0$
and $p>2$, there exists $\delta_0>0$ such that
\[
\|w\|_H\geq\delta_0.
\]
Consequently, every nontrivial profile satisfies
\[
\|w^i\|_H\geq\delta_0.
\]

On the other hand, the weak convergence in each translated
frame yields the successive norm splittings
\[
\|v_n^{i-1}\|_H^2
=
\|w^i\|_H^2+\|v_n^i\|_H^2+o(1).
\]
Since $u_n\rightharpoonup u^0$ in $H$, we also have
\[
\|u_n\|_H^2
=
\|u^0\|_H^2+\|v_n^0\|_H^2+o(1).
\]
Therefore, for every fixed $k\in\mathbb N$, summing the
successive norm splittings gives
\[
\|u_n\|_H^2
=
\|u^0\|_H^2
+\sum_{i=1}^k\|w^i\|_H^2
+\|v_n^k\|_H^2+o(1).
\]
The boundedness of $\{u_n\}$ and the uniform lower bound
for the norms of the nontrivial profiles imply that the
number of profiles is finite. Indeed, if the extraction
process could continue indefinitely, then, for every
$k\in\mathbb N$, letting $n\to\infty$ in the preceding
identity would give
\[
k\delta_0^2
\leq
\limsup_{n\to\infty}\|u_n\|_H^2,
\]
which is impossible for arbitrarily large $k$.

Consequently, after finitely many extractions and passing
to a subsequence if necessary, the vanishing alternative
must occur for some finite $\ell\geq0$, that is,
\[
v_n^\ell\to0\quad\text{strongly in }H.
\]
If $\ell=0$, this is precisely the first alternative of
the splitting lemma.
Thus, the extraction process terminates after finitely
many steps, say at some $\ell\in\mathbb N$. By the
stopping criterion, the final remainder satisfies
\[
v_n^\ell\to0\qquad\text{strongly in }H.
\]
Otherwise, Lemma A.4 would yield, after passing to a
subsequence, another nontrivial critical profile, which
contradicts the finite bound on the number of profiles.
Hence
\[
v_n^\ell\to0\qquad\text{strongly in }H.
\]

By the definition of the remainders, we obtain
\[
u_n
=
u^0+\sum_{i=1}^{\ell}w^i(\cdot-y_n^i)+o(1).
\]
The center escape and pairwise separation properties have
already been established during the extraction process.
Moreover, letting $n\to\infty$ in the successive norm
splittings, we obtain
\[
\|u_n\|_H^2
=
\|u^0\|_H^2
+\sum_{i=1}^{\ell}\|w^i\|_H^2+o(1).
\]

Finally, iterating the energy splitting established in
Lemma A.2, we have
\[
J(u_n)
=
J(u^0)+\sum_{i=1}^{\ell}J(w^i)
+J(v_n^\ell)+o(1).
\]
Since $v_n^\ell\to0$ strongly in $H$ and $J$ is continuous,
we have $J(v_n^\ell)\to J(0)=0$. Therefore,
\[
J(u_n)
=
J(u^0)+\sum_{i=1}^{\ell}J(w^i)+o(1).
\]
Recalling that $J(u_n)\to d$, we conclude that
\[
d=J(u^0)+\sum_{i=1}^{\ell}J(w^i).
\]
This completes the proof.
\end{proof}

\section*{Declaration of Competing Interest}

The authors declare that they have no known competing
financial interests or personal relationships that could
have appeared to influence the work reported in this paper.

\section*{Acknowledgements}\addcontentsline{toc}{section}{Acknowledgements}
This work is supported by the Natural Science Foundation of Chongqing, China (No. CSTB2024NSCQ-MSX0843), the National Natural Science Foundation of China (No. 12361024) and the Team Building Project for Graduate Tutors in Chongqing (No. yds223010).
\phantomsection
\end{document}